\documentclass{scrartcl}
\usepackage{amsmath}
\usepackage{amsfonts}
\usepackage{amssymb}
\usepackage{graphicx}
\usepackage{theorem}
\usepackage{color}%
\usepackage{todonotes}
\usepackage[colorlinks=true]{hyperref}
\usepackage[left=2.5cm,right=2.5cm,vmargin=2.5cm,footnotesep=0.5cm]{geometry}

\newtheorem{thm}{Theorem}[section]
\newtheorem{lm}[thm]{Lemma}
\newtheorem{pr}[thm]{Proposition}
\newtheorem{df}[thm]{Definition}
\newtheorem{rmk}[thm]{Remark}

{\theorembodyfont{\upshape}

}

\let\origitemize\itemize
\def\itemize{\origitemize\itemsep0pt}
\numberwithin{equation}{section}

\makeatletter
\DeclareOldFontCommand{\rm}{\normalfont\rmfamily}{\mathrm}
\DeclareOldFontCommand{\sf}{\normalfont\sffamily}{\mathsf}
\DeclareOldFontCommand{\tt}{\normalfont\ttfamily}{\mathtt}
\DeclareOldFontCommand{\bf}{\normalfont\bfseries}{\mathbf}
\DeclareOldFontCommand{\it}{\normalfont\itshape}{\mathit}
\DeclareOldFontCommand{\sl}{\normalfont\slshape}{\@nomath\sl}
\DeclareOldFontCommand{\sc}{\normalfont\scshape}{\@nomath\sc}
\makeatother
\addtokomafont{disposition}{\rmfamily}

\begin{document}

\title{
Necessary conditions for approximate solutions of vector and set optimization problems with variable domination structure
}

\author{
Marius Durea\thanks{Faculty of Mathematics, Alexandru Ioan
Cuza\ University, 700506--Ia\c{s}i, Romania and Octav Mayer\ Institute of
Mathematics, Ia\c{s}i Branch of Romanian Academy, 700505--Ia\c{s}i, Romania{; e-mail: \texttt{durea@uaic.ro}}}
    \and
    Christian Günther\thanks{Leibniz Universit\"at Hannover, Institut f\"ur Angewandte Mathematik, Welfengarten 1, 30167 Hannover, Germany,
e-mail: \texttt{c.guenther@ifam.uni-hannover.de}, ORCID: 0000-0002-1491-4896}
	\and  Radu Strugariu\thanks{Department of Mathematics, Gheorghe Asachi
Technical University, \ {700506--Ia\c{s}i, Romania and Octav
Mayer\ Institute of Mathematics, Ia\c{s}i Branch of Romanian Academy,
700505--Ia\c{s}i, Romania; e-mail: \texttt{rstrugariu@tuiasi.ro}}}
	\and Christiane Tammer\thanks{Martin Luther University Halle-Wittenberg, Faculty of Natural Sciences II, 
Institute of Mathematics, 06099 Halle (Saale), Germany, e-mail: \texttt{christiane.tammer@mathematik.uni-halle.de}}}

\maketitle

\begin{center}\vspace{-1cm}
	\textbf{Abstract}
\end{center}

\begin{abstract}
	We consider vector and set optimization problems with respect to variable domination structures given by set-valued mappings acting between the preimage space and the image space of the objective mapping, as well as by set-valued mappings with the same input and output space, that coincides with the image space of the objective mapping. The aim of this paper is to derive necessary conditions for approximately nondominated points of problems with a single-valued objective function, employing an extension of Ekeland's Variational Principle for problems with respect to variable domination structures in terms of generalized differentiation in the sense of Mordukhovich. For set-valued objective mappings, we derive necessary conditions for approximately nondominated points of problems with variable domination structure taking into account the incompatibility between openness and optimality and a directional openness result for the sum of set-valued maps. We describe the necessary conditions for approximately nondominated points of set optimization problems with variable domination structure in terms of the limiting (Mordukhovich) generalized differentiation objects.
\end{abstract}

\begin{flushleft}
\end{flushleft}

\begin{flushleft}
	\textbf{Mathematics subject classifications (MSC 2010):} 49J53, 90C29, 92G99
\end{flushleft}

\section{Introduction}

Optimization problems with variable domination structures (vds, for short)
offer a broadly generalized framework for vector optimization (with a fixed
order), which, in turn, extends the scope of classical scalar optimization. It
is worth mentioning that optimization problems with vds have rapidly developed
over the last 15 years, driven by both potential applications and mathematical
interest. The interested reader is referred to \cite{best17}, \cite{BM2014}, \cite{BMST2021}, \cite{CheHuaYan05}, \cite{CheYan02},  \cite{CheYanYu05}, \cite{DST2015}, \cite{E2011}, \cite{Eich2014}, \cite{Engau08}, \cite{EP2016}, 
\cite{yu73dom}, \cite{yu74}, \cite{ZhouTangZhao2024}. A look at this literature shows that, on
the one hand, there are two main ways of considering a vds. Namely, a domination structure is given by a set-valued mapping that acts between the same spaces as the objective mapping, or else, by a set-valued mapping with the same input and output space that coincides with the image space of the objective mapping. On the other
hand, one can consider as objective mapping a single or a set-valued application.

The aim of this paper is to explore some possibilities to obtain {  necessary conditions} for some types of approximate solutions for optimization problems with variable domination structures in the different approaches we mentioned above. To achieve this goal, we establish several new findings in variational and nonlinear analysis, such as variational principles, generalized forms of directional openness for sums of set-valued mappings, and subdifferential calculus.

Therefore, in this work we adapt several principles of optimization in a nontrivial manner to suit the settings under consideration, aiming to address approximate optimality rather than exact optimality.  The mappings involved are inherently nonsmooth, which places the discussion within the framework of nonlinear variational analysis hence the primary tools for expressing optimality conditions are objects of generalized differentiation (see \cite{Mor2006}, \cite{M2018}).

Our investigation is divided into two main parts. The first part focuses on single-valued mappings as objectives. In this case, the two approaches for handling variable domination structures are equivalent, allowing us to center our formalism on just one approach. A key construction, which subsequently enables the derivation of optimality conditions, is a generalized form of the Gerstewitz scalarization functional. On this basis, we present a variant of Ekeland's Variational Principle to problems with respect to (wrt as shortcut) a variable domination structure that acts between the same spaces as the single-valued objective  mapping and then we employ this to derive the necessary conditions for approximate solutions of problems wrt variable domination structures. Extensions of Ekeland's Variational Principle to problems with respect to a variable domination structure for which the preimage as well as the image space are given by the image space of a single-valued objective function are derived in \cite{best17}, \cite{BMST2021} and of a set-valued objective  mapping in \cite{ZhouTangZhao2024}. However, the new scalarization functional we consider here is tailored specifically for our study and introduces additional challenges due to its highly general nature, making it irreducible to known cases.  Nonetheless, by combining the new formulation of the Ekeland Variational Principle with the capabilities of Fréchet--Mordukhovich generalized calculus, we establish {   necessary conditions for approximate solutions}. We also briefly examine a convex setting, where more precise statements about the subdifferential of the scalarization functional can be made.

The second main part addresses the full generality of the proposed problem, specifically considering the case where the objective is a set-valued mapping. In this context, we rely on a well-established optimization principle, which states that a mapping cannot remain open at an efficient point. Building on this, we construct a comprehensive openness assertion for a sum of mappings. This assertion is directional in nature, seeks to provide quantitative estimates for the radii of the balls involved in the conclusion, and incorporates an injectivity condition expressed via the Fréchet coderivative of the set-valued mappings in question. Moreover, its proof uses the classical Ekeland's Variational Principle.  Using this foundation and the aforementioned principle, we derive necessary optimality conditions for problems within the first framework modeled by variable domination structures. Achieving this requires a detailed examination of the constants involved and leveraging the benefits of a finite-dimensional setting. Finally, we propose a method to extend these results to the second framework of vds. This transfer is nontrivial, as several technical assumptions do not directly align with the reduction mechanism between the two approaches to vds.

Our paper is organized as follows: In Section \ref{s-sol_concepts}, we introduce different concepts for approximate solutions of problems wrt variable domination structure and discuss relationships between them. Furthermore, we recall some notions of generalized differentiation in the sense of Mordukhovich. A variational principle of Ekeland's type for problems wrt variable domination structure given by a set-valued  mapping acting between the same spaces as the single-valued objective  mapping and corresponding necessary conditions for approximate solutions are derived in Section \ref{s-single}.  In Section \ref{s-set}, we show a directional openness result for the sum of set-valued  mappings that we will use to derive necessary conditions for approximately nondominated solutions of set-valued problems wrt variable domination structures via certain incompatibility assertions. Finally, we give some conclusions in Section \ref{s-concl}.
\medskip

The notation we use is fairly standard. Therefore, $X,Y$ are normed vector spaces over the real field $\mathbb{R}$, while the topological dual of $X$ is denoted by $X^{\ast}$. We denote by $B(x,r)$ and $\overline{B}(x,r)$ the open and the closed balls of center $x\in X$ and radius $r>0,$ respectively. Also, we denote by $S_X$ the unit sphere of $X$. For a nonempty set $A\subset X,$ we denote by $d\left(  \cdot,A\right)  $ the associated distance function.

\section{Approximate solutions in optimization with vds}\label{s-sol_concepts}

Consider the following set-valued optimization problem:
\[
(P)\hspace{0.4cm}\text{minimize }F(x),\text{ subject to }x\in\Omega,
\]
where $X,Y$ are normed vector spaces, $F:X\rightrightarrows Y$ is a set-valued
map, and $\Omega\subseteq X$ is a closed set. If $C\subseteq Y$ is a proper
closed convex and pointed cone, then it induces a partial order $\leq_{C}$ on
$Y$ by the equivalence $y_{1}\leq_{C}y_{2}$ if and only if $y_{2}-y_{1}\in C.$
The optimality in the usual Pareto sense reads as follows.

\begin{df}
[Pareto solution]\label{df_Pareto}A point $(\overline{x},\overline{y}%
)\in\operatorname{Gr}F\cap\left(  \Omega\times Y\right)  $ is called a local
Pareto solution for $F,$ or for problem $(P)$, if there exists a neighborhood
$U$ of $\overline{x}$ such that%
\[
\left(  F\left(  U\cap\Omega\right)  -\overline{y}\right)  \cap\left(
-C\right)  =\left\{  0\right\}  ,
\]
that is, $\overline{y}$ is a minimal point of $F\left(  U\cap\Omega\right)  $
with respect to $\leq_{C}.$
\end{df}

Now, there are two main ways of generalizing this framework by considering
variable order structures, that is, introducing a domination schema based on a
set-valued mapping. These two approaches differ exactly on the spaces
where such new variable structure is built. More precisely, on the one hand,
one can choose to define the vds on the same spaces as the objective mapping
$F,$ and we call this domination structure $K,$ that is, $K$ is a set-valued
 mapping acting from $X\ $to $Y$ and having as values closed, convex, proper and
pointed cones. On the other hand, one can define the vds as set-valued mapping from $Y$ to $Y,$ with
the same types of values, and we call this domination structure $Q.$

Therefore, from now on, $K:X\rightrightarrows Y$ and $Q:Y\rightrightarrows Y$
are the set-valued  mappings described above, and on this basis one can define
multiple types of solutions that we discuss below.

\begin{df}
[solutions wrt $K$ and $Q$]\label{df_sol}Let $F:X\rightrightarrows Y$ be a
set-valued {map}, $\Omega\subseteq X$ be a closed set, and $(\overline
{x},\overline{y})\in\operatorname*{Gr}F\cap\left(  \Omega\times Y\right)  .$

(i) One says that $\left(  \overline{x},\overline{y}\right)  $ is a local
nondominated solution wrt $K$ for $F$ on {$\Omega$ if} there exists
$\varepsilon>0$ such that, for every $x\in B\left(  \overline{x}%
,\varepsilon\right)  \cap\Omega,$%
\begin{equation}
\left(  F(x)-\overline{y}\right)  \cap\left(  -K(x)\right)  \subseteq\left\{
0\right\}  . \label{nondom1}%
\end{equation}

(ii) One says that $\left(  \overline{x},\overline{y}\right)  $ is a{ local
efficient solution }wrt $K$ for $F$ on $\Omega$ if there exists $\varepsilon
>0$ such that
\[
(F(B\left(  \overline{x},\varepsilon\right)  \cap\Omega)-\overline{y}%
)\cap(-K(\overline{x}))\subseteq\{0\}.
\]

(iii) One says that $\left(  \overline{x},\overline{y}\right)  $ is a local
nondominated solution wrt $Q$ for $F$ on {$\Omega$ if} there exists
$\varepsilon>0$ such that, for every $x\in B\left(  \overline{x}%
,\varepsilon\right)  \cap\Omega,$%
\begin{equation}
\forall y\in F\left(  x\right)  \setminus\left\{  \overline{y}\right\}
:\overline{y}\notin y+Q\left(  y\right)  . \label{nondom2}%
\end{equation}

(iv) One says that $\left(  \overline{x},\overline{y}\right)  $ is a{ local
efficient solution }wrt $Q$ for $F$ on $\Omega$ if there exists $\varepsilon
>0$ such that
\[
(F(B\left(  \overline{x},\varepsilon\right)  \cap\Omega)-\overline{y}%
)\cap(-Q(\overline{y}))\subseteq\{0\}{.}%
\]

If $\Omega:=X,$ we say that $\left(  \overline{x},\overline{y}\right)  $ is a
local nondominated/efficient point for $F$ with respect to $K$ or $Q,$
respectively, and if $\varepsilon=+\infty$ (with the usual convention
$B\left(  \overline{x},+\infty\right)  :=X$), we obtain the global variants of
the above definitions.
\end{df}

{  The concept of domination structures and corresponding solution concepts were introduced by Yu in \cite{yu73dom}, where the sets $Q(y)$ are supposed to be cones for all $y\in Y$. Yu defined a domination structure as a family of cones $Q(y)$, whereas Engau \cite{Engau08} considered it as a set-valued mapping $Q:Y\rightrightarrows Y$. }

\begin{rmk}
\label{rmk_sol}(i) Recently, in \cite{DST2024}, it was shown that approach based
on $Q$ vds can be reduced to that based on $K$ by an extension of the
domain spaces. This reduction will operate as well on the main concepts we are
going to study and will be in force later in this work.

(ii) It is important for the sequel of this work to mention that in the case
where $F$ is a single-valued map, a situation we designate by $F:=f,$ where
$f:X\rightarrow Y$ is a function, we omit $\overline{y}$ and we call the point
$\overline{x}$ a corresponding solution. Moreover, in this case, since for a
fixed $x\in B\left(  \overline{x},\varepsilon\right)  \cap\Omega$ corresponds
exactly one $y=f\left(  x\right)  \in Y,$ there is no difference between the
corresponding $K$ and $Q$ types of solutions.

(iii) If $K(\cdot)\equiv C$ or $Q\left(  \cdot\right)  \equiv C$, where $C\subseteq
Y$ a fixed proper closed convex and pointed cone, then all these solution
concepts coincide with Pareto efficiency from Definition \ref{df_Pareto}.

(iv) In the framework of vds, the main focus is on nondominated solutions,
since the efficient solutions are collapsing to the case of Pareto solutions,
because the cone in the right-hand side of the defining relations is actually
fixed ($K\left(  \overline{x}\right)  $ and $Q\left(  \overline{y}\right)$, respectively).
\end{rmk}

As mentioned, our concern in this work is to study necessary optimality
conditions for approximate solution in optimization problems with vds, and in
view of Remark \ref{rmk_sol} (iv), we define below only nondominated
approximate solutions, that is, the main concepts under current investigation.

\begin{df}
[approximately nondominated solutions wrt $K$ and $Q$]\label{d-appr-solK}Let
$F:X\rightrightarrows Y$ be a set-valued {map}, $\Omega\subseteq X$ be a
closed set, and $(\overline{x},\overline{y})\in\operatorname*{Gr}F\cap\left(
\Omega\times Y\right)  .$ Consider $k\in Y\setminus\{0\},$ and $\varepsilon
,\delta>0$.

(i) One says that\textbf{ }$\left(  \overline{x},\overline{y}\right)  $ is an
$\left(  \varepsilon,\delta,{{k}}\right)  ${-nondominated solution} wrt $K$
for $F$ on $\Omega$ if for every $x\in
B\left(  \overline{x},\varepsilon\right)  \cap\Omega,${%
\begin{equation}
\left(  F\left(  x\right)  -\overline{y}+\delta k\right)  \cap\left(
-K(x)\right)  \subseteq\{0\}. \label{edk-nondom}%
\end{equation}
}

(ii) One says that\textbf{ }$\left(  \overline{x},\overline{y}\right)  $ is an
$\left(  \varepsilon,\delta,{{k}}\right)  ${-nondominated solution} wrt $Q$
for $F$ on $\Omega$ if for every $x\in
B\left(  \overline{x},\varepsilon\right)  \cap\Omega,${%
\begin{equation}
\forall y\in F\left(  x\right)  \setminus\left\{  \overline{y}\right\}
:\;\overline{y}-\delta k\not \in y+(Q(y)\setminus\{0\}).\label{edk-nondom_Q}%
\end{equation}
Similar observations as above, concerning global, unrestricted solutions or
the case where }$F:=f${, are in order.}
\end{df}

\begin{rmk}
If one takes $\delta:=0$ in the above definition, the concepts of
{nondominated solutions} wrt $K$ and $Q$ are obtained, respectively. Most of
the results stated in the sequel also hold in this case, with obvious
modifications. In this work, we concentrate only on the approximate case.
\end{rmk}

\begin{rmk}
\label{reduction}An $\left(  \varepsilon,\delta,{{k}}\right)  $-nondominated
solution wrt $Q$ can be reduced to an $\left(  \varepsilon,\delta,{{k}%
}\right)  $-nondominated solution wrt $K$ by the following device (see
\cite{DST2024}). Suppose, in the above notation, that $\left(  \overline
{x},\overline{y}\right)  $ an $\left(  \varepsilon,\delta,{{k}}\right)
$-nondominated solution wrt $Q$ for $F$ on $\Omega.$ Consider $\widetilde
{F},\widetilde{K}:X\times Y\rightrightarrows Y$ the set-valued  mappings given by%
\begin{align}
\widetilde{F}\left(  x,y\right)   &  :=\left\{
\begin{array}
[c]{ll}%
\left\{  y\right\}   & \text{if }y\in F\left(  x\right)  \\
\emptyset & \text{otherwise,}%
\end{array}
\right.  \label{F_tilt}\\
\widetilde{K}\left(  x,y\right)   &  :=Q\left(  y\right)  \text{ for all }y\in
Y,\label{K_tilt}%
\end{align}
and remark that $\left(  \left(  \overline{x},\overline{y}\right)
,\overline{y}\right)  \in\operatorname*{Gr}\widetilde{F}$ is an $\left(
\left(  \varepsilon,+\infty\right)  ,\delta,{{k}}\right)  $-nondominated
solution wrt $\widetilde{K}$ for $\widetilde{F}$ on $\Omega\times Y$ if for
every $\left(  x,y\right)  \in\left(  B\left(  \overline{x},\varepsilon
\right)  \cap\Omega\right)  \times Y,$%
\[
\left(  \widetilde{F}(x,y)-\overline{y}+\delta k\right)  \cap\left(
-\widetilde{K}(x,y)\right)  \subseteq\left\{  0\right\}  ,
\]
and this reduces to: for any $x\in B\left(  \overline{x},\varepsilon\right)
\cap\Omega,$%
\[
\forall y\in F\left(  x\right)  :\left(  y-\overline{y}+\delta k\right)
\cap\left(  -Q\left(  y\right)  \right)  \subseteq\left\{  0\right\}  ,
\]
which is actually (\ref{edk-nondom_Q}).
\end{rmk}

\begin{rmk}
\label{rmk_app_f}When $F:=f$, the {nondominated solutions} wrt $K$ and
$Q\ $actually coincide (see also Remark \ref{rmk_sol} (ii)).
\end{rmk}

The {  necessary conditions for approximate solutions} we are going to present are expressed in terms of
generalized differentiation objects in the sense of Mordukhovich (see,
\cite{Mor2006}, \cite{M2018}). A very short description of these objects which
are essential in our study is given next.

For a nonempty subset $\Omega\ $of the Asplund space $X$ and $x\in\Omega,$ the
Fr\'{e}chet normal cone to $\Omega$ at $x$ is%

\begin{equation}
\widehat{N}(\Omega,x):=\left\{  x^{\ast}\in X^{\ast}\mid\underset
{u\overset{\Omega}{\rightarrow}x}{\lim\sup}\frac{\langle x^{\ast},u-x\rangle}{\left\Vert
u-x\right\Vert }\leq0\right\}  , \label{Fr_cone}%
\end{equation}
where $u\overset{\Omega}{\rightarrow}x$ means that $u\rightarrow x$ and
$u\in\Omega.$ If $x\notin\Omega,$ we let $\widehat{N}(\Omega,x):=\emptyset.$
If $\Omega$ is closed around $\overline{x}$, the limiting (Mordukhovich)
normal cone is given by%
\begin{equation}
N(\Omega,\overline{x}):=\left\{  x^{\ast}\in X^{\ast}\mid\exists x_{n}%
\overset{\Omega}{\rightarrow}\overline{x},x_{n}^{\ast}\overset{\ast
}{\rightarrow}x^{\ast},x_{n}^{\ast}\in\widehat{N}(\Omega,x_{n}),\forall
n\in\mathbb{N}\right\}  , \label{M_cone}%
\end{equation}
where $\overset{\ast}{\rightarrow}$ means the convergence in the weak$^{\ast}$ topology.

Next, we present the associated coderivative constructions for set-valued
{maps}. Let $F:X\rightrightarrows Y$ be a set-valued {map} with the domain and
the graph defined by%
\[
\operatorname{Dom}F:=\{x\in X\mid F(x)\neq\emptyset\}\quad\text{and}%
\quad\operatorname{Gr}F:=\{(x,y)\mid y\in F(x)\},
\]
and $(\overline{x},\overline{y})\in\operatorname*{Gr}F.$ Then the Fr\'{e}chet
coderivative at $(\overline{x},\overline{y})$ is the set-valued {map}
$\widehat{D}^{\ast}F(\overline{x},\overline{y}):Y^{\ast}\rightrightarrows
X^{\ast}$ given by
\begin{equation}
\widehat{D}^{\ast}F(\overline{x},\overline{y})(y^{\ast}):=\{x^{\ast}\in
X^{\ast}\mid(x^{\ast},-y^{\ast})\in\widehat{N}(\operatorname{Gr}%
F,(\overline{x},\overline{y}))\}, \label{Fr_coder}%
\end{equation}
while the Mordukhovich coderivative of $F$ at $(\overline{x},\overline{y})$ is
the set-valued {map} $D^{\ast}F(\overline{x},\overline{y}):Y^{\ast
}\rightrightarrows X^{\ast}$ given by
\begin{equation}
D^{\ast}F(\overline{x},\overline{y})(y^{\ast}):=\{x^{\ast}\in X^{\ast}%
\mid(x^{\ast},-y^{\ast})\in N(\operatorname{Gr}F,(\overline{x},\overline
{y}))\}. \label{M_coder}%
\end{equation}
As usual, when $F:=f$ is a function, since $\overline{y}\in F\left(
\overline{x}\right)  $ means $\overline{y}=f\left(  \overline{x}\right)  ,$ we
write $\widehat{D}^{\ast}f\left(  \overline{x}\right)  $ for $\widehat
{D}^{\ast}f\left(  \overline{x},\overline{y}\right)  ,$ and similarly for
$D^{\ast}$.

For a function $f:X\rightarrow\mathbb{R\cup\{+\infty\}}$ finite at
$\overline{x}\in X\ $and lower semicontinuous around $\overline{x},$ the
Fr\'{e}chet subdifferential of $f$ at $\overline{x}$ is the set
\[
\widehat{\partial}f(\overline{x}):=\{x^{\ast}\in X^{\ast}\mid(x^{\ast}%
,-1)\in\widehat{N}(\operatorname*{epi}f,(\overline{x},f(\overline{x})))\},
\]
where $\operatorname*{epi}f$ denotes the epigraph of $f$, while the limiting
(Mordukhovich) subdifferential of $f$ at $\overline{x}$ is given, according to
\cite{Mor2006}, by%
\[
\partial f(\overline{x}):=\{x^{\ast}\in X^{\ast}\mid(x^{\ast},-1)\in
N(\operatorname*{epi}f,(\overline{x},f(\overline{x})))\}.
\]
It is well-known that if $f$ is a convex function, then $\widehat{\partial
}f(\overline{x})$ and $\partial f(\overline{x})$ coincide with the Fenchel
subdifferential. However, in general, $\widehat{\partial}f(\overline
{x})\subseteq\partial f(\overline{x}),$ and the following generalized Fermat
rule holds: if $\overline{x}\in X$ is a local minimum point for
$f:X\rightarrow\mathbb{R\cup\{+\infty\}}$, then $0\in\widehat{\partial
}f(\overline{x}).$ Some other properties of these objects will be reminded and
used when needed.

\section{Optimality conditions for single-valued objective mappings}\label{s-single}

In this section, we consider the particular case $F:=f$ (in the above mentioned
notation) and we get {  necessary conditions for approximately} {nondominated}
solutions with respect to $K$ (see also Remark \ref{rmk_app_f}). 

{  Roughly speaking, in an Ekeland-type variational principle, assuming boundedness from below and  lower semicontinuity of the objective mapping, one has to show that for an approximate solution of the original problem there exists another approximate solution (condition (i)) such that the new approximate solution belongs to a neighborhood of the initial approximate solution (condition (ii)), and is an exact solution of a slightly perturbed problem (condition (iii)). Condition (iii) is useful for deriving necessary conditions for approximate solutions. The assertions in a variational principle of Ekeland's type are shown without any convexity or compactness assumptions. In the next theorem, the variational principle for $\left(  \varepsilon,\delta,{{k}}\right)  ${-nondominated solution} wrt $K:X\rightrightarrows Y$ (see Definition \ref{d-appr-solK})
for a single-valued mapping $f:X\rightarrow Y$ on $\Omega$ is shown under weaker assumptions in comparison with corresponding results in the literature, where so far there are only variational principles for problems with variable dominance structures $Q:Y\rightrightarrows Y$, see \cite[Theorem 4.7]{best17} and \cite[Theorem 3.9]{BMST2021}.}

The main
device will be an adaptation of the classical Ekeland Variational Principle
(EVP, for short) to an appropriate scalar functional.

More precisely, we employ the following scalarizing function (as an extended
version of the Gerstewitz scalarizing function): $s_{k,a-K}^{f}:X\rightarrow
\mathbb{R}\cup\{\pm\infty\}$ {defined} for $a\in Y,$ $k\in Y\setminus\{0\},$
$f:X\rightarrow Y$ and $K:X\rightrightarrows Y$ by the formula
\begin{equation}
s_{k,a-K}^{f}(x):=\operatorname*{inf}\left\{  t\in\mathbb{R}\mid f(x)\in
a+tk-K(x)\right\}  .\label{def:sK}%
\end{equation}

Here, we use the convention $\inf\emptyset:=+\infty.$ Notice that in our
setting, where $K\left(  x\right)  $ is a closed set, the above $\inf$ is
actually a minimum if it is not $\pm\infty$. Recall that the Gerstewitz
scalarizing function {   is defined for a {proper, closed} subset $R$ of a linear topological space $Y$ and a nonzero vector $k \in Y$ with $k \in 0^+R$, where 
$0^+ R:=\{u\in Y\mid \forall y\in R,\ \forall t\in
\mathbb{R}_{+}:\; y+tu\in R\}$
is the recession cone of the nonempty set $R\subset Y$, as }
$s_{k,R}%
:Y\rightarrow\mathbb{R}\cup\{\pm\infty\}$,
\begin{equation}\label{f-skR}
s_{k,R}(y):=\operatorname*{inf}%
\left\{  t\in\mathbb{R}\mid y\in tk-R\right\}  .
\end{equation}
This functional, along the
newly introduced one, will intervene in our discussion later.

Now, we formulate in the next theorem an assertion for domination structures
$K:X\rightrightarrows Y$, which is an adaptation of the EVP in this situation
via the functional defined by (\ref{def:sK}).

\begin{thm}
\label{thm:bestK} Let $X$ be a Banach space, while $Y$ is a normed vector
space, $k\in Y\setminus\{0\}$, let $K:X\rightrightarrows Y$ be a vds,
$\Omega\subseteq X$ be a closed set, and let $f:X\rightarrow Y$ be a
vector-valued mapping. Given $\varepsilon,\delta>{0}$, consider an $\left(
\varepsilon,\delta,{{k}}\right)  $-nondominated solution $\tilde{x}$ with
respect to $K$ for $f$ on $\Omega.$ Impose the following assumptions:

(a) \textsc{(Continuity condition)} $f$ is continuous over $\overline
{B}\left(  \tilde{x},\varepsilon\right)  \cap\Omega$, and the domination
mapping $K$ has locally closed graph;

(b) \textsc{(scalarization condition)} $K(x)+(0,+\infty)k\subseteq
K(x)\setminus\left\{  0\right\}  $ for all $x\in\overline{B}\left(
\tilde{x},\varepsilon\right)  \cap\Omega$.

Then for all $\varepsilon^{\prime}\in\left(  0,\varepsilon\right)  $ there
exists an element $\overline{x}\in\Omega$ (depending on $\varepsilon$ and
$\delta$) such that the following conditions hold:

(i) $s_{k,f\left(  \tilde{x}\right)  -K}^{f}(\overline{x})+\sqrt{{\delta}%
}\left\Vert \overline{x}-\tilde{x}\right\Vert \leq s_{k,f\left(  \tilde
{x}\right)  -K}^{f}(\tilde{x})$.

(ii) $\left\Vert \overline{x}-\tilde{x}\right\Vert \leq\min\left\{
\sqrt{\delta},\varepsilon^{\prime}\right\}  $.

(iii) $\overline{x}$ is a strict minimal point on $\overline{B}\left(
\tilde{x},\varepsilon^{\prime}\right)  \cap\Omega$ of the scalar function
$s_{k,f\left(  \tilde{x}\right)  -K}^{f}\left(  \cdot\right)  +\sqrt{\delta
}\left\Vert \overline{x}-\cdot\right\Vert $.
\end{thm}

\noindent\textbf{Proof.} As mentioned, we employ the functional given in
(\ref{def:sK}), taking $a=f\left(  \tilde{x}\right)  $ and showing that the
assumptions of the classical EVP are fulfilled for this function on the metric
space $\overline{B}\left(  \tilde{x},\varepsilon^{\prime}\right)  \cap\Omega$,
where $\varepsilon^{\prime}\in\left(  0,\varepsilon\right)  .$

First, we show that this function is bounded from below by $-\delta$ on
$\overline{B}\left(  \tilde{x},\varepsilon^{\prime}\right)  \cap\Omega.$ Assume, by way
of contradiction, the existence of $x\in\overline{B}\left(  \tilde
{x},\varepsilon^{\prime}\right)  \cap\Omega$ such that $s_{k,f\left(
\tilde{x}\right)  -K}^{f}(x)<-\delta.$ Then there is $t\in\mathbb{R}$,
$t<-\delta$ such that
\begin{equation}
f(\tilde{x})+tk-f(x)\in K(x).\label{f-f1}%
\end{equation}
Since, by the scalarization condition {(b),} $K(x)+(-t-\delta)k\subseteq
K(x)\setminus\{0\}$, we have
\[
f(\tilde{x})-\delta k\in f(x)+K(x)+(-t-\delta)k\subseteq f(x)+(K(x)\setminus
\{0\}),
\]
in contradiction to the $\left(  \varepsilon,\delta,{{k}}\right)
$-nondominatedness of $\tilde{x}$. So, indeed,
\begin{equation}
\forall
x\in\overline{B}\left(  \tilde{x},\varepsilon^{\prime}\right)  \cap
\Omega: s_{k,f\left(  \tilde{x}\right)  -K}^{f}(x)\geq-\delta.\label{f-s1}%
\end{equation}
On the other hand, $s_{k,f\left(  \tilde{x}\right)  -K}^{f}(\tilde{x})=0,$
that is, $\operatorname*{inf}\left\{  t\in\mathbb{R}\mid tk\in K(\tilde
{x})\right\}  =0$ since $K(\tilde{x})$ is a cone, and the {scalarization
condition easily yields }that $(-\infty,0)k\cap K(\tilde{x})=\emptyset$
(otherwise, it would exist $\alpha<0$ such that $\alpha k\in K(\tilde{x}),$
whence $0\in-\alpha k+K(\tilde{x})\subseteq K(\tilde{x})\setminus\left\{
0\right\}  $){.}

Therefore,%

\[
0=s_{k,f\left(  \tilde{x}\right)  -K}^{f}(\tilde{x})\leq\inf_{x\in \overline{B}\left(
\tilde{x},\varepsilon^{\prime}\right)  \cap\Omega}s_{k,f\left(  \tilde
{x}\right)  -K}^{f}({x})+\delta,
\]
i.e., $\tilde{x}$ is a $\delta$-minimal solution of $s_{k,f\left(  \tilde
{x}\right)  -K}^{f}$ over $\overline{B}\left(  \tilde{x},\varepsilon^{\prime
}\right)  \cap\Omega$.

The next step is to show that the functional $s_{k,f\left(  \tilde{x}\right)
-K}^{f}$ is lower semicontinuous, by showing the closedness of its level sets
with respect to the metric space $\overline{B}\left(  \tilde{x},\varepsilon
^{\prime}\right)  \cap\Omega$. Take $t\in\mathbb{R}$ and consider the
corresponding level set of $s_{k,f\left(  \tilde{x}\right)  -K}^{f},$ that
is,
\[
A_{t}:=\left\{  x\in\overline{B}\left(  \tilde{x},\varepsilon^{\prime}\right)
\cap\Omega\mid s_{k,f\left(  \tilde{x}\right)  -K}^{f}(x)\leq t\right\}  .
\]
If $A_{t}=\emptyset,$ there is nothing to prove. Otherwise, take any
convergent sequence $\left(  x_{n}\right)  $ in $A_{t}$ and denote by
$x_{\ast}$ its limit. For every nonzero natural number $n$, there is $t_{n}%
\in\left(  -2\delta,t+n^{-1}\right)  $ such that
\[
f\left(  \tilde{x}\right)  +t_{n}k-f(x_{n})\in K(x_{n}).
\]
Since $\left(  t_{n}\right)  $ is a bounded sequence of real numbers, one can
suppose, without loss of generality, that it is convergent to a number
$\tilde{t}\leq t.$ The continuity condition (a) ensures that $f(x_{n}%
)\rightarrow f\left(  x_{\ast}\right)  $ and $f\left(  \tilde{x}\right)
+\tilde{t}k-f(x_{\ast})\in K(x_{\ast}),$ so $x_{\ast}\in A_{t},$ whence the
latter set is closed.

Consequently, we have shown that all assumptions of the classical EVP are
fulfilled for $s_{k,f\left(  \tilde{x}\right)  -K}^{f}$ over $\overline
{B}\left(  \tilde{x},\varepsilon^{\prime}\right)  \cap\Omega$. From this
result we obtain the existence of $\overline{x}\in\overline{B}\left(
\tilde{x},\varepsilon^{\prime}\right)  \cap\Omega$ such that%
\begin{align*}
&  s_{k,f\left(  \tilde{x}\right)  -K}^{f}(\overline{x})+\sqrt{{\delta}%
}\left\Vert \tilde{x}-\overline{x}\right\Vert \leq s_{k,f\left(  \tilde
{x}\right)  -K}^{f}(\tilde{x}),\\
&  \left\Vert \tilde{x}-\overline{x}\right\Vert \leq\sqrt{\delta},\\
&  s_{k,f\left(  \tilde{x}\right)  -K}^{f}(x)+\sqrt{\delta}\left\Vert
x-\overline{x}\right\Vert >s_{k,f\left(  \tilde{x}\right)  -K}^{f}%
(\overline{x}) \mbox{ for all } \; x\in\overline{B}\left(  \tilde{x}%
,\varepsilon^{\prime}\right)  \cap\Omega\setminus\left\{  \overline
{x}\right\} .
\end{align*}
This completes the proof.\hfill$\square$

\begin{rmk}
The scalarization condition is clearly satisfied if
\[
k\in%
{\displaystyle\bigcap\limits_{x\in \overline{B}\left(  \tilde{x},\varepsilon\right)
\cap\Omega}}
K(x)\setminus\left\{  0\right\}  .
\]

\end{rmk}

{   \begin{rmk}
Even if the statement (iii) in Theorem \ref{thm:bestK} is shown in a scalarized form, it can be used to derive necessary conditions for approximate solutions.    
\end{rmk}}
Now we are able to derive necessary conditions for $\left(  \varepsilon
,\delta,{{k}}\right)  $-nondominated solutions in the framework assumed in
this section by employing Theorem \ref{thm:bestK}.

\begin{thm}
\label{nc-effK} Suppose that $X$ and $Y$ are Asplund spaces and assume the
notation and the hypotheses of Theorem \ref{thm:bestK}. Let $\overline{x}%
\in\Omega$ be the element from the conclusion of Theorem \ref{thm:bestK}. Then
the following statements hold:

(i) if $\sqrt{\delta}\geq\varepsilon$ and $\overline{x}-\tilde{x}%
\notin-N\left(  \Omega,\overline{x}\right)  $, then $0\in\partial s_{k,f\left(
\tilde{x}\right)  -K}^{f}({\overline{x}})+N(\Omega,\overline{x})+\sqrt{\delta
}\overline{B}_{X^{\ast}}+\operatorname*{cone}\left\{  \overline{x}-\tilde
{x}\right\}  ;$

(ii) if $\sqrt{\delta}<\varepsilon$, then $0\in\partial s_{k,f\left(  \tilde
{x}\right)  -K}^{f}({\overline{x}})+N(\Omega,\overline{x})+\sqrt{\delta
}\overline{B}_{X^{\ast}}.$
\end{thm}

\noindent\textbf{Proof. }In both cases we employ the statement of Theorem
\ref{thm:bestK} and the usual generalized calculus for the limiting (Mordukhovich) subdifferential (see \cite[Chapter 3, especially Theorem 3.36]{Mor2006}), to write
\[
0\in\partial\left(  s_{k,f\left(  \tilde{x}\right)  -K}^{f}\left(
\cdot\right)  +\sqrt{\delta}\left\Vert \overline{x}-\cdot\right\Vert \right)
\left(  \overline{x}\right)  +N\left(  \overline{B}\left(  \tilde
{x},\varepsilon^{\prime}\right)  \cap\Omega,\overline{x}\right)  .
\]
The difference between the two cases is that in the second one, the inequality
$\sqrt{\delta}<\varepsilon$ means that one can find $\varepsilon^{\prime}%
\in\left(  0,\varepsilon\right)  $ such that $\sqrt{\delta}<\varepsilon
^{\prime},$ so, according to the second conclusion of Theorem \ref{thm:bestK},
$\overline{x}\ $is an interior point of $\overline{B}\left(  \tilde
{x},\varepsilon^{\prime}\right)  .$ Therefore, in this case $N\left(
\overline{B}\left(  \tilde{x},\varepsilon^{\prime}\right)  \cap\Omega
,\overline{x}\right)  =N\left(  \Omega,\overline{x}\right)  $. The inequality
of the first case allows $\overline{x}$ to be on the border of $\overline
{B}\left(  \tilde{x},\varepsilon^{\prime}\right)  ,$ so in this situation,
under the additional assumption $\overline{x}-\tilde{x}\notin-N\left(
\Omega,\overline{x}\right)  ,$ one has $N\left(  \overline{B}\left(  \tilde
{x},\varepsilon^{\prime}\right)  \cap\Omega,\overline{x}\right)  \subseteq
N\left(  \Omega,\overline{x}\right)  +\operatorname*{cone}\left\{
\overline{x}-\tilde{x}\right\}  $ (see \cite[Theorem 1.21, Corollary
3.5]{Mor2006}). The conclusion follows.\hfill$\square$

{  \begin{rmk}
A necessary condition for approximately nondominated solutions of vector optimization problems wrt variable domination structures $Q:Y\rightrightarrows Y$ in terms of the limiting (Mordukhovich) subdifferential is derived in \cite[Theorem 5.5]{best17} under different assumptions than in Theorem \ref{nc-effK}.    
\end{rmk}}

As one can see in Theorem \ref{nc-effK}, it is important to have an estimation
of the subdifferential of the function $s_{k,a-K}^{f}$ defined by relation
(\ref{def:sK}). Of course, in every specific case, an ad-hoc computation can
be envisaged, but it is important to mention as well a general situation when
something more precise can be said. As usual, this situation is a
(generalized) convex setting. The standing assumption we make here is that
$k\in%
{\displaystyle\bigcap\limits_{x\in X}}
K(x)\setminus\left\{  0\right\}  ,$ which, in particular, means that the set 
$P:=
{\displaystyle\bigcap\limits_{x\in X}}
K(x)$ is a proper closed convex and pointed cone.

\begin{df}
Let $G:X\rightrightarrows Y$ be a set valued  mapping and $R\subseteq X$ be a
proper closed convex and pointed cone. Then one says that $G$ is convex wrt
$R$ if for all $\alpha\in\left(  0,1\right)  $ and $x_{1},x_{2}\in X$, one has
\[
\alpha G\left(  x_{1}\right)  +\left(  1-\alpha\right)  G\left(  x_{2}\right)
\subseteq G\left(  \alpha x_{1}+\left(  1-\alpha\right)  x_{2}\right)  +R.
\]

\end{df}

Of course, this definition applies with obvious changes for $G:=g,$ a
single-valued map.

\begin{pr}
Suppose that $f$ and $K$ are convex wrt $P$ and $s_{k,a-K}^{f}$ is proper (this
happens, for instance, if $k\in\operatorname*{int}P$). Then

(i) $s_{k,a-K}^{f}$ is convex;

(ii) for all $\overline{x}\in X$ and $x^{\ast}\in\partial s_{k,f\left(
\overline{x}\right)  -K}^{f}\left(  \overline{x}\right)  ,$ there is $y^{\ast
}\in P^{+}$ with $y^{\ast}\left(  k\right)  =1$ and
\[
x^{\ast}\in\partial\left(  y^{\ast}\circ\left(  f\left(  \cdot\right)
-f\left(  \overline{x}\right)  \right)  \right)  \left(  \overline{x}\right)
.
\]

\end{pr}

\noindent\textbf{Proof. }(i) Let $\alpha\in\left(  0,1\right)  $ and
$x_{1},x_{2}\in X.$ If $s_{k,a-K}^{f}\left(  x_{1}\right)  =+\infty$ or
$s_{k,a-K}^{f}\left(  x_{2}\right)  =+\infty$, there is nothing to prove.
Otherwise,
\begin{align*}
f\left(  x_{1}\right)   &  \in a+s_{k,a-K}^{f}\left(  x_{1}\right)  k-K\left(
x_{1}\right)  ,\\
f\left(  x_{2}\right)   &  \in a+s_{k,a-K}^{f}\left(  x_{2}\right)  k-K\left(
x_{2}\right)  ,
\end{align*}
whence
\begin{align*}
f\left(  \alpha x_{1}+\left(  1-\alpha\right)  x_{2}\right)   &  \in\alpha
f\left(  x_{1}\right)  +\left(  1-\alpha\right)  f\left(  x_{2}\right)  -P\\
&  \subseteq a+\left(  \alpha s_{k,a-K}^{f}\left(  x_{1}\right)  +\left(
1-\alpha\right)  s_{k,a-K}^{f}\left(  x_{2}\right)  \right)  k-\left(  \alpha
K\left(  x_{1}\right)  +\left(  1-\alpha\right)  K\left(  x_{2}\right)
\right)  -P\\
&  \subseteq a+\left(  \alpha s_{k,a-K}^{f}\left(  x_{1}\right)  +\left(
1-\alpha\right)  s_{k,a-K}^{f}\left(  x_{2}\right)  \right)  k-K\left(  \alpha
x_{1}+\left(  1-\alpha\right)  x_{2}\right)  -P-P\\
&  \subseteq a+\left(  \alpha s_{k,a-K}^{f}\left(  x_{1}\right)  +\left(
1-\alpha\right)  s_{k,a-K}^{f}\left(  x_{2}\right)  \right)  k-K\left(  \alpha
x_{1}+\left(  1-\alpha\right)  x_{2}\right)  .
\end{align*}
This shows that $s_{k,a-K}^{f}\left(  \alpha x_{1}+\left(  1-\alpha\right)
x_{2}\right)  \leq\alpha s_{k,a-K}^{f}\left(  x_{1}\right)  +\left(
1-\alpha\right)  s_{k,a-K}^{f}\left(  x_{2}\right)  ,$ that is, the convexity
of $s_{k,a-K}^{f}.$

(ii) Take now $\overline{x}\in X$ and $x^{\ast}\in\partial s_{k,f\left(
\overline{x}\right)  -K}^{f}\left(  \overline{x}\right)  ,$ meaning that
\[ \forall x\in X:  
s_{k,f\left(  \overline{x}\right)  -K}^{f}\left(  x\right)  -s_{k,f\left(
\overline{x}\right)  -K}^{f}\left(  \overline{x}\right)  \geq\left\langle
x^{\ast},x-\overline{x}\right\rangle.
\]
In particular, since $P\subseteq K\left(  x\right)  $, one has $s_{k,f\left(
\overline{x}\right)  -P}^{f}\left(  x\right)  \geq s_{k,f\left(  \overline
{x}\right)  -K}^{f}\left(  x\right)  $ (where we see $P$ as a constant
set-valued  mapping from $X$ to $Y$), therefore
\begin{equation}\label{f-subineq}
\forall x\in X: s_{k,f\left(  \overline{x}\right)  -P}^{f}\left(  x\right)  -0\geq\left\langle
x^{\ast},x-\overline{x}\right\rangle .
\end{equation}
Clearly, the former function depends on $x$ only via $f,$ so we can see it as
the composition $s_{k,P}\circ t_{-f\left(  \overline{x}\right)  }\circ f,$
where $t_{-f\left(  \overline{x}\right)  }$ is the translation $y\mapsto
y-f\left(  \overline{x}\right)  $, and the functional $s_{k,P}$ is given by \eqref{f-skR} for $R:= P$. Inequality \eqref{f-subineq} is actually
equivalent to $x^{\ast}\in\partial s_{k,f\left(  \overline{x}\right)  -P}%
^{f}\left(  \overline{x}\right)  .$ So, $x^{\ast}\in\partial\left(
s_{k,P}\circ t_{-f\left(  \overline{x}\right)  }\circ f\right)  \left(
\overline{x}\right)  $, and we are now using \cite[Theorem 2.8.10]{Zal} to
write
\[
x^{\ast}\in%
{\displaystyle\bigcup}
\left\{  \partial\left(  y^{\ast}\circ t_{-f\left(  \overline{x}\right)
}\circ f\right)  \left(  \overline{x}\right)  \mid y^{\ast}\in P^{+}%
\cap\partial s_{k,P}\left(  0\right)  \right\}  .
\]
Employing the form of $\partial s_{k,P}\left(  0\right)  $ obtained in
\cite{DT2009}, we get that for $x^{\ast}$, there is $y^{\ast}\in P^{+}$ with
$y^{\ast}\left(  k\right)  =1$ such that $x^{\ast}\in\partial\left(  y^{\ast
}\circ\left(  f\left(  \cdot\right)  -f\left(  \overline{x}\right)  \right)
\right)  \left(  \overline{x}\right)  $. This finishes the proof.\hfill
$\square$

\section{Optimality conditions for set-valued objective mappings}\label{s-set}

\subsection{Necessary conditions for {  approximately nondominated} points wrt $K$}

This section is devoted to   necessary conditions for   approximately
nondominated solutions in full generality, which means that we tackle the case
of set-valued objective mappings. The main ingredients are an instance of the
well-known incompatibility between openness and optimality and a
directional openness result for the sum of a finite number of set-valued maps.

\bigskip

In the next proposition, we use the indicator set-valued  mapping $\Delta_{\Omega
}:X\rightrightarrows Y$, given by%
\[
\Delta_{\Omega}\left(  x\right)  :=\left\{
\begin{array}
[c]{ll}%
\left\{  0\right\}   & \text{if }x\in\Omega\\
\emptyset & \text{otherwise.}%
\end{array}
\right.
\]
The next incompatibility result holds for approximate minima.

\begin{pr}
\label{prop_incom}Let $\varepsilon\in\left(  0,\infty\right]  ,\delta
\in\left(  0,\infty\right)  $, and $k\in\left(
{\displaystyle\bigcap_{x\in B\left(  \overline{x},\varepsilon\right)
\cap\Omega}}
K\left(  x\right)  \right)  \setminus\left\{  0\right\}  .$ If $\left(
\overline{x},\overline{y}\right)  \in\operatorname*{Gr}F\cap\left(
\Omega\times Y\right)  $ is an $\left(  \varepsilon,\delta,{{k}}\right)
${-nondominated solution }wrt $K$ for $F$ on $\Omega,$ then one cannot have
$\delta^{\prime}>\delta$ such that%
\begin{equation}
B\left(  \overline{y},\delta^{\prime}\right)  \cap\left[  \overline
{y}-\operatorname*{cone}\left\{  k\right\}  \right]  \subseteq\left(
F+K+\Delta_{\Omega}\right)  \left(  B\left(  \overline{x},\varepsilon\right)
\right)  . \label{dir_op}%
\end{equation}

\end{pr}

\noindent\textbf{Proof.} Observe that, since $K\left(  x\right)  $ is a convex
cone, relation (\ref{edk-nondom}) is equivalent to%
\[
\forall x\in B\left(  \overline{x},\varepsilon\right)  \cap\Omega :\;\left(
F\left(  x\right)  +K\left(  x\right)  -\overline{y}+\delta k\right)
\cap\left(  -K(x)\right)  \subseteq\left\{  0\right\}  .
\]
Suppose, by contradiction, that there is $\delta^{\prime}>\delta$ such that
(\ref{dir_op}) holds. This implies, if we take arbitrary $\delta^{\prime
\prime}\in\left(  \delta,\delta^{\prime}\right)  ,$ that%
\[
\overline{y}-\delta^{\prime\prime}k\in B\left(  \overline{y},\delta^{\prime
}\right)  \cap\left[  \overline{y}-\operatorname*{cone}\left\{  k\right\}
\right]  \in\left(  F+K+\Delta_{\Omega}\right)  \left(  B\left(  \overline
{x},\varepsilon\right)  \right)  ,
\]
so there exists $x_{0}\in B\left(  \overline{x},\varepsilon\right)  \cap
\Omega$ such that
\[
-\delta^{\prime\prime}k\in F\left(  x_{0}\right)  +K\left(  x_{0}\right)
-\overline{y},
\]
hence%
\[
-\left(  \delta^{\prime\prime}-\delta\right)  k\in F\left(  x_{0}\right)
+K\left(  x_{0}\right)  +\delta k-\overline{y}.
\]
Since obviously $-\left(  \delta^{\prime\prime}-\delta\right)  k\in-K\left(
x_{0}\right)  \setminus\left\{  0\right\}  ,$ we obtain a
contradiction.$\hfill\square$

\bigskip

Using the previous proposition, we see that it is of interest the formulation
of sufficient conditions for the directional openness of the sum
$F+K+\Delta_{\Omega}.$

First of all, let us precisely define the concept of directional openness. For
supplemental details, the reader is referred to \cite{DPS2017}.

\begin{df}
Let $u\in Y \setminus\{0\}$. One says that $F$ is directionally linearly open at $(\overline{x}%
,\overline{y})\in\operatorname*{Gr}F$ with respect to $\left\{  u\right\}  $
with modulus $\alpha>0$ if there is $\varepsilon>0$ such that, for every
$r\in(0,\varepsilon),$%
\begin{equation}
B(\overline{y},\alpha r)\cap\lbrack\overline{y}-\operatorname*{cone}\left\{
u\right\}  ]\subseteq F(B(\overline{x},r)). \label{linop_dir}%
\end{equation}

\end{df}

In the proof of the main openness result, we will use the minimal time
function, which we present next. Let $A\subseteq X$ be a nonempty set and
$u\in S_{X}$. Then the minimal time function, defined and systematically
analyzed in \cite{Nam-Zal-2013}, is given as follows: for any $x\in X,$%
\[
T_{u}(x,A):=\inf\left\{  t\geq0\mid x+tu\in A\right\}.
\]

The following properties of the minimal time function will be important in the
sequel (see, e.g. \cite{Nam-Zal-2013}, \cite{DPS2016}, \cite{DPS2017}). 

\begin{pr}
\label{pr_basic}(i) The domain of the directional minimal time function with
respect to $u$ is%
\[
\operatorname*{dom}{T_{u}(\cdot,A)=A-\operatorname*{cone}}\left\{  u\right\}
.
\]

(ii) For any $x\in X,$%
\begin{equation}
d(x,A)\leq T_{u}(x,A).\label{ineq_d_t}%
\end{equation}

(iii) If $A,B\subseteq X,$ then%
\begin{equation}
T_{u}(x,A\cup B)=\min\left\{  T_{u}(x,A),T_{u}(x,B)\right\}  . \label{prb4}%
\end{equation}

(iv) For any $A_{1},A_{2}\subseteq X,$ the following relation holds:
\[
T_{u}(x_{1}+x_{2},A_{1}+A_{2})\leq T_{u}(x_{1},A_{1})+T_{u}(x_{2},A_{2}%
),\quad\forall\ x_{i}\in A_{i}-\operatorname{cone}\left\{  u\right\}
,\ i=1,2.
\]

(v) If $A$ is closed, or if $X$ is reflexive and $A$ is weakly closed, then:
\[
\forall x\in A-{\operatorname*{cone}}\left\{  u\right\}  :x+T_{u}(x,A)u\in A.
\]

(vi) Let $\lambda>0$. Then:%
\[
{[T_{u}(\cdot,A)<\lambda]=A-[0,\lambda)\cdot}\left\{  u\right\}
\]
and
\[
{A-[0,\lambda]\cdot\left\{  u\right\}  \subseteq[T_{u}(\cdot,A)\leq
\lambda]\subseteq%
{\displaystyle\bigcap_{\varepsilon>0}}
\left(  A-[0,\lambda+\varepsilon)\cdot\left\{  u\right\}  \right)  }.
\]

Moreover, if $A$ is closed, then
\[
\lbrack T_{u}(\cdot,A)\leq\lambda]=A-[0,\lambda]\cdot\left\{  u\right\}  ,
\]
and, in particular, $T_{u}$ is lower semicontinuous.

(vii) If $A$ is a convex set, then $T_{u}(\cdot,A)$ is a convex function.

(viii) Suppose $A$ is closed and take $\overline{x}\in(A-{\operatorname*{cone}%
}\left\{  u\right\}  )\setminus A.$ Then for any $a\in A$ with $\overline
{x}+T_{u}(\overline{x},A)u=a,$ one has
\begin{equation}
\widehat{\partial}T_{u}(\cdot,A)(\overline{x})\subseteq\{x^{\ast}\in X^{\ast
}\mid\left\langle x^{\ast},u\right\rangle =-1\}\cap\widehat{N}(A,a).
\label{subdiff_tmin}%
\end{equation}

\end{pr}

In the case $A:=\left\{  v\right\}  ,$ for some $v\in X,$ we denote
$T_{u}(x,\left\{  v\right\}  )$ by $T_{u}(x,v).$

\bigskip

Another ingredient refers to the alliedness of sets, which is presented next (for more details and historical comments, see
\cite{DHNS2014}).

\begin{df}
\label{df_alliedness}If $S_{1},...,S_{p}$ are subsets of a normed vector space
$X,$ closed around $\overline{x}\in S_{1}\cap...\cap S_{p}$, one says that
they are allied at $\overline{x}$ (for the Fr\'{e}chet normal cones) whenever
$(x_{in})\overset{S_{i}}{\rightarrow}\overline{x},x_{in}^{\ast}\in\widehat
{N}(S_{i},x_{in}),i=\overline{1,p},$ the relation $(x_{1n}^{\ast}%
+...+x_{pn}^{\ast})\rightarrow0$ implies $(x_{in}^{\ast})\rightarrow0$ for
every $i=\overline{1,p}$ .
\end{df}

Let us formulate now the directional openness result for the sum of set-valued
maps, on which we will build later in formulating optimality conditions via the
incompatibility assertion from Proposition \ref{prop_incom}.

\begin{thm}
\label{suf_cond_dir_op_sum} Suppose that $X,Y$ are finite dimensional spaces,
$F_{i}:X\rightrightarrows Y,i=\overline{1,p}$ are closed-graph multifunctions,
$(\overline{x},\overline{y}_{i})\in\operatorname*{Gr}F_{i},i=\overline{1,p},$
and $u\in S_{Y}.$ Assume that the sets given by%
\begin{equation}
D_{i}:=\{(x,y_{1},...,y_{p})\in X\times Y^{p}\mid y_{i}\in F_{i}(t)\},i=\overline
{1,p}\label{D1p}%
\end{equation}
are allied at $(\overline{x},\overline{y}_{1},...,\overline{y}_{p})$, and there
exists $c>0$ and $r>0$ such that for every ${(u_{i},v_{i})\in}%
\operatorname*{Gr}F_{i}\cap\lbrack B(\overline{x},r)\times B(\overline{y}%
_{i},r)],$ ${i=}\overline{{1,p}},$ every $y^{\ast}\in Y^{\ast},$ with
$\left\langle y^{\ast},u\right\rangle =1,z_{i}^{\ast}\in2cB_{Y^{\ast}%
},i=\overline{{2,p}}$,
 and every $x_{1}^{\ast}\in\widehat{D}^{\ast}F_{1}%
(u_{1},v_{1})(y^{\ast}),$ $x_{i}^{\ast}\in\widehat{D}^{\ast}F_{i}(u_{i}%
,v_{i})(y^{\ast}-z_{i}^{\ast}),$ $i=\overline{{2,p}},$%
\begin{equation}
c\leq\Vert x_{1}^{\ast}+...+x_{p}^{\ast}\Vert.\label{ineq_x12}%
\end{equation}

Then for every $a\in(0,c),$ $H:X\rightrightarrows Y$ given as%
\begin{equation}
H\left(  x\right)  :=%
{\displaystyle\sum\limits_{i=1}^{p}}
F_{i}\left(  x\right)  \label{H_sum}%
\end{equation}
is directionally open at $(\overline{x},\overline{y}_{1}+...+\overline{y}%
_{p})\in\operatorname*{Gr}H$ with respect to $\left\{  u\right\}  $ with
modulus $a.$
\end{thm}

\noindent\textbf{Proof.} Fix $a\in\left(  0,c\right)  ,$ and denote $D:=%
{\displaystyle\bigcap\limits_{i=1}^{p}}
D_{i},\overline{y}:=\overline{y}_{1}+...+\overline{y}_{p}.$ We will find
$\theta>0$ such that, for every $\rho\in(0,\theta)$,%
\begin{equation}
B(\overline{y},a\rho)\cap\left[  \overline{y}-\operatorname*{cone}\left\{
u\right\}  \right]  \subseteq H(B(\overline{x},\rho)).\label{interm}%
\end{equation}
Choose $b\in(0,1)$ such that $\dfrac{a}{a+1}<b<\dfrac{c}{c+1}\ $. Then one can find $\theta>0$ such that the following are satisfied:

\begin{itemize}
\item $b^{-1}a\theta<r;$

\item $\theta<r;$

\item $\dfrac{a}{a+1}<b+2\theta<\dfrac{c}{c+1}.$
\end{itemize}

Fix $\rho\in(0,\theta)$. Choose now $v\in B(\overline{y},\rho a)\cap\left[
\overline{y}-\operatorname*{cone}\left\{  u\right\}  \right]  ,$ hence
$T_{u}\left(  v,\overline{y}\right)  <\infty.$ Observe that the set $D$ is
closed because of the closedness of the graphs of $F_{1},...,F_{p}.$ Endow the
space $X\times Y^{p}$ with the distance%
\[
d((x,y_{1},...,y_{p}),(x^{\prime},y_{1}^{\prime},...,y_{p}^{\prime
})):=b\left(  \left\Vert x-x^{\prime}\right\Vert +\left\Vert y_{1}%
-y_{1}^{\prime}\right\Vert +...+\left\Vert y_{p}-y_{p}^{\prime}\right\Vert
\right)
\]
and apply the Ekeland Variational Principle for the lower semicontinuous and
lower bounded function $f:D\rightarrow\mathbb{R},$
\[
f(x,y_{1},y_{2},...,y_{p}):=T_{u}(v,y_{1}+y_{2}+...+y_{p})=T_{-u}(y_{1}%
+y_{2}+...+y_{p},v)
\]
at $(\overline{x},\overline{y}_{1},...,\overline{y}_{p})\in\operatorname*{dom}%
f.$\ Denote $y:=y_{1}+y_{2}+...+y_{p}.$ We find $(x_{0},y_{10},...,y_{p0})\in
D$ such that%
\begin{equation*}
f(x_{0},y_{10},...,y_{p0})\leq f(\overline{x},\overline{y}_{1},...,\overline
{y}_{p})-d((x_{0},y_{10},...,y_{p0}),(\overline{x},\overline{y}_{1}%
,...,\overline{y}_{p})), \label{ek1b}%
\end{equation*}
and%
\begin{equation*}
f(x_{0},y_{10},...,y_{p0})\leq f(x,y_{1},...,y_{p})-d((x_{0},y_{10}%
,...,y_{p0}),(x,y_{1},...,y_{p})),\quad\forall(x,y_{1},...,y_{p})\in D.
\label{ek2b}%
\end{equation*}
Equivalently, denoting also $y_{0}:=y_{10}+y_{20}+...+y_{p0},$ this can be
written as%
\begin{equation}
T_{u}(v,y_{0})\leq T_{u}(v,\overline{y})-b\left(  \left\Vert \overline
{x}-x_{0}\right\Vert +\left\Vert \overline{y}_{1}-y_{10}\right\Vert
+...+\left\Vert \overline{y}_{p}-y_{p0}\right\Vert \right)  \label{ek_n1}%
\end{equation}
and%
\begin{equation}
T_{u}(v,y_{0})\leq T_{u}(v,y)+b(\left\Vert x-x_{0}\right\Vert +\left\Vert
y_{1}-y_{10}\right\Vert +...+\left\Vert y_{p}-y_{p0}\right\Vert ),\quad
\forall(x,y_{1},...,y_{p})\in D. \label{ek_n2}%
\end{equation}

Remark that $T_{u}(v,\overline{y})=\left\Vert v-\overline{y}\right\Vert ,$
hence by (\ref{ek_n1}) $T_{u}(v,y_{0})$ is finite and, consequently, $v\in
y_{0}-\operatorname*{cone}\left\{  u\right\}  $ and $T_{u}(v,y_{0})=\left\Vert
v-y_{0}\right\Vert .$ Moreover, relation (\ref{ek_n1}) can be written as%
\begin{equation}
\left\Vert v-y_{0}\right\Vert \leq\left\Vert v-\overline{y}\right\Vert
-b\left(  \left\Vert \overline{x}-x_{0}\right\Vert +\left\Vert \overline
{y}_{1}-y_{10}\right\Vert +...+\left\Vert \overline{y}_{p}-y_{p0}\right\Vert
\right)  .\label{ek_1b}%
\end{equation}
This implies that%
\begin{align*}
\left\Vert \overline{x}-x_{0}\right\Vert  &  \leq b^{-1}\left\Vert
v-\overline{y}\right\Vert <b^{-1}a\rho<b^{-1}a\theta<r,\\
\left\Vert \overline{y}_{1}-y_{10}\right\Vert +...+\left\Vert \overline{y}%
_{p}-y_{p0}\right\Vert  &  \leq b^{-1}\left\Vert v-\overline{y}\right\Vert
<b^{-1}a\rho<b^{-1}a\theta<r,
\end{align*}
hence%
\[
(x_{0},y_{10},...,y_{p0})\in B(\overline{x},r)\times B(\overline{y}%
_{1},r)\times...\times B(\overline{y}_{p},r).
\]

If $v=y_{0},$ then%
\begin{align*}
b\left\Vert \overline{x}-x_{0}\right\Vert  &  \leq\left\Vert v-\overline
{y}\right\Vert -b\left(  \left\Vert \overline{y}_{1}-y_{10}\right\Vert
+...+\left\Vert \overline{y}_{p}-y_{p0}\right\Vert \right)  \\
&  \leq\left\Vert v-\overline{y}\right\Vert -b\left(  \left\Vert \left(
\overline{y}_{1}+...+\overline{y}_{p}\right)  -\left(  y_{10}+...+y_{p0}%
\right)  \right\Vert \right)  \\
&  =\left(  1-b\right)  \left\Vert v-\overline{y}\right\Vert <(1-b)a\rho
<b\rho,
\end{align*}
hence $x_{0}\in B(\overline{x},\rho)$ and $v=y_{0}\in H(x_{0})\subseteq
H(B(\overline{x},\rho)),$ which proves the desired conclusion.

We want to prove now that $v=y_{0}$ is the only possible situation. Suppose
then by contradiction that $v\not =y_{0}$, and consider the functions
$g,h:X\times Y^{p}\rightarrow\mathbb{R},$%
\begin{align*}
g(x,y_{1},...,y_{p}) &  :=T_{u}(v,y_{1}+...+y_{p}),\\
h(x,y_{1},...,y_{p}) &  :=b(\left\Vert x-x_{0}\right\Vert +\left\Vert
y_{1}-y_{10}\right\Vert +...+\left\Vert y_{p}-y_{p0}\right\Vert ).
\end{align*}
From (\ref{ek_n2}), we have that the point $(x_{0},y_{10},...,y_{p0})$ is a
minimum point for $g+h$ on the set $D$, or, equivalently, $(x_{0}%
,y_{10},...,y_{p0})$ is a global minimum point for the function $g+h+\delta
_{D}$, where by $\delta_{D}$ we denote the indicator function of the set $D,$ i.e.,
the function which values $0$ for the elements of $D,$ and $+\infty$ otherwise. Applying the generalized Fermat rule, we have that%
\[
(0,0,...,0,0)\in\widehat{\partial}(g+h+\delta_{D})(x_{0},y_{10},...,y_{p0}).
\]

Using the fact that all the functions $g,h$ and $\delta_{D}$ are lower
semicontinuous, and since $X,Y$ are finite dimensional, we can apply the weak
fuzzy calculus rule for the Fr\'{e}chet subdifferential (see, e.g.,
\cite[Theorem 3]{Ioffe1983}, \cite{Fab1989}).

Choose $\gamma\in(0,\min\{\rho,2^{-1}r\})$ such that $v\notin\overline
{B}(y_{0},p\gamma)=y_{10}+...+y_{p0}+p\overline{B}(0,\gamma)=\overline
{B}(y_{10},\gamma)+...+\overline{B}(y_{p0},\gamma)$ and obtain that there
exist
\begin{align*}
(u_{1},v_{11},...,v_{p1})  &  \in\overline{B}(x_{0},\gamma)\times\overline
{B}(y_{10},\gamma)\times...\times\overline{B}(y_{p0},\gamma)\text{ with }%
v\in\left(  v_{11}+...+v_{p1}\right)  -\operatorname*{cone}\left\{  u\right\}
,\\
(u_{2},v_{12},...v_{p2})  &  \in\overline{B}(x_{0},\gamma)\times\overline
{B}(y_{10},\gamma)\times...\times\overline{B}(y_{p0},\gamma),\\
(u_{3},v_{13},...v_{p3})  &  \in\lbrack\overline{B}(x_{0},\gamma
)\times\overline{B}(y_{10},\gamma)\times...\times\overline{B}(y_{p0}%
,\gamma)]\cap D
\end{align*}
such that%
\[
(0,0,...,0,0)\in\widehat{\partial}g(u_{1},v_{11},...,v_{p1})+\widehat
{\partial}h(u_{2},v_{12},...v_{p2})+\widehat{\partial}\delta_{D}(u_{3}%
,v_{13},...v_{p3})+\gamma\left(  \overline{B}_{X^{\ast}}\times\overline
{B}_{Y^{\ast}}\times...\times\overline{B}_{Y^{\ast}}\right)  .
\]
Observe that $h$ is the sum of $p+1$ convex functions, Lipschitz on $X,$ hence
$\widehat{\partial}h$ coincides with the sum of the convex subdifferentials.
Hence, we obtain from before%
\begin{align*}
(0,0,...,0,0)  &  \in\left\{  0\right\}  \times\widehat{\partial}T_{u}%
(v,\cdot+...+\cdot)\left(  v_{11},...,v_{p1}\right)  +b\overline{B}_{X^{\ast}%
}\times\{0\}\times...\times\{0\}\times\{0\}\\
&  +\{0\}\times b\overline{B}_{Y^{\ast}}\times...\times\{0\}+...+\{0\}\times
\{0\}\times...\times b\overline{B}_{Y^{\ast}}\\
&  +\widehat{N}(D,(u_{3},v_{13},...,v_{p3}))+\rho\overline{B}_{X^{\ast}}%
\times\rho\overline{B}_{Y^{\ast}}\times...\times\rho\overline{B}_{Y^{\ast}}.
\end{align*}
Observe now that the function%
\[
\left(  y_{1},...,y_{p}\right)  \mapsto T_{u}(v,y_{1}+...+y_{p})=T_{-u}%
(y_{1}+...+y_{p},v)
\]
is the composition between the convex function $T_{u}(v,\cdot)$ and the linear
 mapping $S\left(  y_{1},...,y_{p}\right)  :=y_{1}+...+y_{p},$ so the Fr\'{e}chet
subdifferential reduces to the convex subdifferential $\partial.$ We apply the
calculus rule for the convex subdifferential and we obtain that%
\[
\widehat{\partial}T_{-u}(\cdot+...+\cdot,v)\left(  v_{11},...,v_{p1}\right)
=S^{\ast}\left(  \partial T_{-u}(\cdot,v)\left(  v_{11}+...+v_{p1}\right)
\right)  ,
\]
where $S^{\ast}:Y^{\ast}\rightarrow Y^{\ast}\times...\times Y^{\ast}$ is the
adjoint of $S$ and is given by $S^{\ast}\left(  y^{\ast}\right)  =\left(
y^{\ast},...,y^{\ast}\right)  $ for any $y^{\ast}\in Y^{\ast}.$

Also, because $v_{11}+...+v_{p1}\in v-\operatorname*{cone}\left\{  -u\right\}
$ and $v\neq v_{11}+...+v_{p1}\in\overline{B}(y_{0},p\gamma),$ we get using
(\ref{subdiff_tmin}) that
\[
\widehat{\partial}T_{-u}(\cdot,v)(v_{11}+...+v_{p1})\subseteq\{y^{\ast}\in
Y^{\ast}\mid\left\langle y^{\ast},u\right\rangle =1\}.
\]
Hence, the relation from above can be written as%
\begin{align*}
(0,0,...,0,0)  &  \in\left\{  0\right\}  \times\left\{  \left(  y^{\ast
},...,y^{\ast}\right)  \mid\left\langle y^{\ast},u\right\rangle =1\right\}
+\widehat{N}(D,(u_{3},v_{13},...,v_{p3}))\\
&  +(b+\rho)\left(  \overline{B}_{X^{\ast}}\times\overline{B}_{Y^{\ast}}%
\times...\times\overline{B}_{Y^{\ast}}\right)  .
\end{align*}
Now, use the alliedness of $D_{1},...,D_{p}$ at $(\overline{x},\overline
{y}_{1},...,\overline{y}_{p})$ to get that%
\begin{align*}
\widehat{N}(D,(u_{3},v_{13},...,v_{p3}))  &  \subseteq\widehat{N}(D_{1}%
,(u_{4},v_{14},...,v_{p4}))+\widehat{N}(D_{2},(u_{5},v_{15},...,v_{p5}))\\
&  +...+\widehat{N}(D_{p},(u_{p+3},v_{1,p+3},...,v_{p,p+3}))+\rho\left(
\overline{B}_{X^{\ast}}\times\overline{B}_{Y^{\ast}}\times...\times
\overline{B}_{Y^{\ast}}\right)  ,
\end{align*}
where
\[
(u_{i+3},v_{1,i+3},...,v_{p,i+3})\in\lbrack\overline{B}(u_{3},\gamma)\times
\overline{B}(v_{13},\gamma)\times...\times\overline{B}(v_{p3},\gamma)]\cap
D_{i},i=\overline{1,p}.
\]
It follows that%
\begin{align*}
(0,0,...,0,0)  &  \in\left\{  0\right\}  \times\left\{  \left(  y^{\ast
},...,y^{\ast}\right)  \mid\left\langle y^{\ast},u\right\rangle =1\right\}
+\widehat{N}(D_{1},(u_{4},v_{14},...,v_{p4}))\\
&  +\widehat{N}(D_{2},(u_{5},v_{15},...,v_{p5}))+...+\widehat{N}%
(D_{p},(u_{p+3},v_{1,p+3},...,v_{p,p+3}))\\
&  +(b+2\rho)\left(  \overline{B}_{X^{\ast}}\times\overline{B}_{Y^{\ast}%
}\times...\times\overline{B}_{Y^{\ast}}\right)  .
\end{align*}

In conclusion, there exist
\begin{align*}
&  y_{0}^{\ast}\in Y^{\ast}\text{ with }\left\langle y_{0}^{\ast
},u\right\rangle =1,\text{ hence }\left\Vert y_{0}^{\ast}\right\Vert \geq1\\
&  (u_{4}^{\ast},v_{14}^{\ast},0,...,0)\in\widehat{N}(D_{1},(u_{4}%
,v_{14},...,v_{p4}))\Leftrightarrow u_{4}^{\ast}\in\widehat{D}^{\ast}%
F_{1}(u_{4},v_{14})(-v_{14}^{\ast}),\\
&  (u_{5}^{\ast},0,v_{25}^{\ast},...,0)\in\widehat{N}(D_{2},(u_{5}%
,v_{15},...,v_{p5}))\Leftrightarrow u_{5}^{\ast}\in\widehat{D}^{\ast}%
F_{2}(u_{5},v_{25})(-v_{25}^{\ast}),\\
&  ...\\
&  (u_{p+3}^{\ast},0,...,0,v_{p,p+3}^{\ast},0)\in\widehat{N}(D_{p}%
,(u_{p+3},v_{1,p+3},...,v_{p,p+3}))\Leftrightarrow u_{p+3}^{\ast}\in\widehat
{D}^{\ast}F_{p}(u_{p+3},v_{p,p+3})(-v_{p,p+3}^{\ast}),\\
&  (u_{p+4}^{\ast},v_{1,p+4}^{\ast},...,v_{p,p+4}^{\ast})\in\overline
{B}_{X^{\ast}}\times\overline{B}_{Y^{\ast}}\times...\times\overline
{B}_{Y^{\ast}}%
\end{align*}
such that%
\begin{align*}
&  -u_{4}^{\ast}-u_{5}^{\ast}-...-u_{p+3}^{\ast}-(b+2\rho)u_{p+4}^{\ast}=0,\\
&  -y_{0}^{\ast}-v_{14}^{\ast}-(b+2\rho)v_{1,p+4}^{\ast}=0,\\
&  ...\\
&  -y_{0}^{\ast}-v_{p,p+3}^{\ast}-(b+2\rho)v_{p,p+4}^{\ast}=0.
\end{align*}

Observe that%
\[
\left\langle y_{0}^{\ast}+(b+2\rho)v_{1,p+4}^{\ast},u\right\rangle
\geq1-(b+2\rho)>1-(b+2\theta)>0,
\]
and denote%
\begin{align*}
x_{1}^{\ast} &  :=\left\langle y_{0}^{\ast}+(b+2\rho)v_{1,p+4}^{\ast
},u\right\rangle ^{-1}u_{4}^{\ast},\\
&  ...\\
x_{p}^{\ast} &  :=\left\langle y_{0}^{\ast}+(b+2\rho)v_{1,p+4}^{\ast
},u\right\rangle ^{-1}u_{p+3}^{\ast},\\
y^{\ast} &  :=-\left\langle y_{0}^{\ast}+(b+2\rho)v_{1,p+4}^{\ast
},u\right\rangle ^{-1}v_{14}^{\ast}\\
&  =\left\langle y_{0}^{\ast}+(b+2\rho)v_{1,p+4}^{\ast},u\right\rangle
^{-1}\left(  y_{0}^{\ast}+(b+2\rho)v_{1,p+4}^{\ast}\right)  ,\\
z_{2}^{\ast} &  :=y^{\ast}+\left\langle y_{0}^{\ast}+(b+2\rho)v_{1,p+4}^{\ast
},u\right\rangle ^{-1}v_{25}^{\ast}\\
&  =\left\langle y_{0}^{\ast}+(b+2\rho)v_{1,p+4}^{\ast},u\right\rangle
^{-1}(b+2\rho)\left(  v_{1,p+4}^{\ast}-v_{2,p+4}^{\ast}\right)  \\
&  ...\\
z_{p}^{\ast} &  :=y^{\ast}+\left\langle y_{0}^{\ast}+(b+2\rho)v_{1,p+4}^{\ast
},u\right\rangle ^{-1}v_{p,p+3}^{\ast}\\
&  =\left\langle y_{0}^{\ast}+(b+2\rho)v_{1,p+4}^{\ast},u\right\rangle
^{-1}(b+2\rho)  \left(  v_{1,p+4}^{\ast}-v_{p,p+4}^{\ast}\right)
  .
\end{align*}

In conclusion, one has%
\begin{align*}
&  \left\langle y^{\ast},u\right\rangle =1\\
&  \left\Vert z_{i}^{\ast}\right\Vert =(b+2\rho)\frac{\left\Vert
v_{1,p+4}^{\ast}-v_{i,p+4}^{\ast}\right\Vert }{\left\langle y_{0}^{\ast
}+(b+2\rho)v_{1,p+4}^{\ast},u\right\rangle }\leq\frac{2(b+2\rho)}{1-(b+2\rho
)}<2\frac{b+2\theta}{1-(b+2\theta)}<2c,i=\overline{{2,p}}\\
&  x_{1}^{\ast}\in\widehat{D}^{\ast}F_{1}(u_{4},v_{14})(y^{\ast}),\\
&  x_{i}^{\ast}\in\widehat{D}^{\ast}F_{i}(u_{i+3},v_{i,i+3})(y^{\ast}%
-z_{i}^{\ast}),i=\overline{{2,p}}%
\end{align*}
where%
\begin{equation}
\left\Vert x_{1}^{\ast}+...+x_{p}^{\ast}\right\Vert =\frac{\left\Vert
u_{4}^{\ast}+...+u_{p+3}^{\ast}\right\Vert }{\left\langle y_{0}^{\ast
}+(b+2\rho)v_{1,p+4}^{\ast},u\right\rangle }=\frac{\left\Vert (b+2\rho
)u_{p+4}^{\ast}\right\Vert }{\left\langle y_{0}^{\ast}+(b+2\rho)v_{1,p+4}%
^{\ast},u\right\rangle }\leq\frac{b+2\rho}{1-(b+2\rho)}<c.\label{e.2}%
\end{equation}

Remark that $(u_{4},v_{14})\in\operatorname*{Gr}F_{1}$ and
\begin{align*}
(u_{4},v_{14})  &  \in\overline{B}(u_{4},\gamma)\times\overline{B}%
(v_{14},\gamma)\subseteq\overline{B}(x_{0},2\gamma)\times\overline{B}%
(y_{10},2\gamma)\\
&  \subseteq B(\overline{x},r)\times B(\overline{y}_{1},r).
\end{align*}
Similarly, $(u_{i+3},v_{i,i+3})\in\operatorname*{Gr}F_{i}\cap\lbrack
B(\overline{x},r)\times B(\overline{y}_{i},r)],i=\overline{{2,p}}.$ Using now
(\ref{ineq_x12}) and (\ref{e.2}), we get that%
\[
c\leq\left\Vert x_{1}^{\ast}+...+x_{p}^{\ast}\right\Vert <c,
\]
a contradiction. It follows that $v=z_{0}$ is the only possible situation and
the conclusion follows.$\hfill\square$

\bigskip

The next lemma is proven in \cite{DST2024} in a more general form.

\begin{lm}
\label{alliedness}Let $X,Y$ be finite dimensional spaces, $\Omega\subseteq X$
be a closed set, $F,K:X\rightrightarrows Y$ be closed-graph set-valued mappings,
and $\left(  \overline{x},\overline{y}\right)  \in X\times Y$ such that
$\overline{x}\in\Omega,(\overline{x},\overline{y})\in\operatorname*{Gr}%
F,\left(  \overline{x},0\right)  \in\operatorname*{Gr}K.$ Moreover, {suppose}
that the following assumptions are satisfied:%
\begin{align*}
&  D^{\ast}K(\overline{x},0)(0)=\left\{  0\right\}  ,\\
&  D^{\ast}F(\overline{x},\overline{y})(0)\cap\left(  -N\left(  \Omega
,\overline{x}\right)  \right)  =\left\{  0\right\}  .
\end{align*}
Then the sets:%
\begin{align}
D_{1} &  :=\{(x,y,z,w)\in X\times Y^{3}:y\in F(x)\}\nonumber\\
D_{2} &  :=\{(x,y,z,w)\in X\times Y^{3}:z\in K(x)\}\label{D13}\\
D_{3} &  :=\{(x,y,z,w)\in X\times Y^{3}:w\in\Delta_{\Omega}\left(  x\right)
\}\nonumber
\end{align}
are allied at $\left(  \overline{x},\overline{y},0,0\right)  .$
\end{lm}

The following result, given in \cite{DST2015}, will be useful in obtaining the
positivity of a multiplier when formulating optimality conditions.

\begin{lm}
\label{pos}Let $K:X\rightrightarrows Y$ be a vds, and $\left(  \overline
{x},\overline{y}\right)  \in\operatorname*{Gr}K.$ If $\widehat{D}^{\ast
}K(\overline{x},\overline{y})(y^{\ast})\neq\emptyset,$ then $y^{\ast}\in
K\left(  \overline{x}\right)  ^{+}.$
\end{lm}

We are in the position to formulate the main openness result that
subsequently will allow to formulate optimality conditions.

\begin{thm}
\label{dir_op_FKD}Let $X,Y$ be finite dimensional spaces, $\Omega\subseteq X$
be a closed set, $F:X\rightrightarrows Y$ be a closed-graph set-valued mapping,
$K:X\rightrightarrows Y$ be a closed-graph vds, and $\left(  \overline
{x},\overline{y}\right)  \in X\times Y$ such that $\overline{x}\in
\Omega,(\overline{x},\overline{y})\in\operatorname*{Gr}F,\left(  \overline
{x},0\right)  \in\operatorname*{Gr}K.$ Denote $P:=%
{\displaystyle\bigcap_{x\in B\left(  \overline{x},\eta\right)  }}
K\left(  x\right)  $ for some $\eta>0$, and suppose $k\in P\cap S_{Y}.$
Moreover, suppose that the following conditions are satisfied:

(a) (\textsc{Transversality Conditions})
\begin{align*}
&  D^{\ast}K(\overline{x},0)(0)=\left\{  0\right\}  ,\\
&  D^{\ast}F(\overline{x},\overline{y})(0)\cap\left(  -N\left(  \Omega
,\overline{x}\right)  \right)  =\left\{  0\right\}  .
\end{align*}

(b) (\textsc{Injectivity Condition}) there exist $r\in\left(  0,\eta\right)
,c>0$ such that for every $y^{\ast}\in P^{+}$ with $\left\langle y^{\ast
},k\right\rangle =1,$ every $(x_{1},y_{1})\in\operatorname*{Gr}F\cap\lbrack
B(\overline{x},r)\times B(\overline{y},r)],$ $(x_{2},y_{2})\in
\operatorname*{Gr}K\cap\lbrack B(\overline{x},r)\times B(0,r)],$ $x_{3}%
\in\Omega\cap B\left(  \overline{x},r\right)  ,$ every $z^{\ast}%
\in2cB_{Y^{\ast}},$ and every $x_{1}^{\ast}\in\widehat{D}^{\ast}F(x_{1}%
,y_{1})(y^{\ast}-z^{\ast}),x_{2}^{\ast}\in\widehat{D}^{\ast}K(x_{2}%
,y_{2})(y^{\ast})$ and $x_{3}^{\ast}\in\widehat{N}\left(  \Omega
,x_{3}\right)  ,$ one has%
\[
c\leq\left\Vert x_{1}^{\ast}+x_{2}^{\ast}+x_{3}^{\ast}\right\Vert .
\]

\noindent Take arbitrary $a\in(0,c),$ and $\theta$ such that
\begin{equation}
0<\theta<\min\left\{  r,\frac{1}{4}\left(  \dfrac{c}{c+1}-\dfrac{a}%
{a+1}\right)  ,\frac{r}{2a}\left(  \dfrac{c}{c+1}+\dfrac{a}{a+1}\right)
\right\}  . \label{epsilon}%
\end{equation}
Then for any $\rho\in\left(  0,\theta\right)  ,$%
\[
B(\overline{y},a\rho)\cap\left[  \overline{y}-\operatorname*{cone}\left\{
k\right\}  \right]  \subseteq\left(  F+K+\Delta_{\Omega}\right)  \left(
B(\overline{x},\rho)\right)  .
\]

\end{thm}

\noindent\textbf{Proof.} Using Lemma \ref{alliedness}, we obtain that the sets
given by (\ref{D13}) are allied at $\left(  \overline{x},\overline
{y},0,0\right)  .$

Let us discuss why it is sufficient to take here $y^{\ast}\in P^{+}.$ If
$\widehat{D}^{\ast}K(x_{2},y_{2})(y^{\ast})=\emptyset,$ then the
(\textsc{Injectivity Condition}) trivially holds. If $\widehat{D}^{\ast
}K(x_{2},y_{2})(y^{\ast})\neq\emptyset,$ then by Lemma \ref{pos} we have that $y^{\ast}\in
K\left(  x_{2}\right)  ^{+}.$ But, since $x_{2}\in B(\overline{x},r)\subseteq
B\left(  \overline{x},\eta\right)  ,$ we have that $P\subseteq K\left(
x_{2}\right)  ,$ hence by the antimonotonicity of the polar, it follows that
$y^{\ast}\in P^{+}.$

Next, by choosing%
\[
b:=\frac{1}{2}\left(  \dfrac{c}{c+1}+\dfrac{a}{a+1}\right)  \in\left(
\dfrac{a}{a+1},\dfrac{c}{c+1}\right)  ,
\]
we observe that, for any $\theta$ given by (\ref{epsilon}), the inequalities:

\begin{itemize}
\item $b^{-1}a\theta<r;$

\item $\theta<r;$

\item $\dfrac{a}{a+1}<b+2\theta<\dfrac{c}{c+1}$
\end{itemize}

\noindent are satisfied. The conclusion follows from Theorem
\ref{suf_cond_dir_op_sum} applied for $F_{1}:=K,$ $F_{2}:=F$ and
$F_{3}:=\Delta_{\Omega}.$ This choice is the appropriate one, since we want to
ensure $y^{\ast}\in P^{+}.\hfill\square$

\bigskip

Let us discuss now the above quantitative openness result from the point of
view of the radii of the balls involved in the incompatibility result stated
in Proposition \ref{prop_incom}.

Take $c=r$ in the above result (since one can choose $\min$ $\left\{
r,c\right\}  $ as the common value). Denote $m_{a}$ the minimum from the
right-hand side of formula (\ref{epsilon}), that is,
\[
m_{a}:=\min\left\{  c,\frac{1}{4}\left(  \dfrac{c}{c+1}-\dfrac{a}{a+1}\right)
,\frac{c}{2a}\left(  \dfrac{c}{c+1}+\dfrac{a}{a+1}\right)  \right\}.
\]
The question of interest at this point is the following one: fix
$\varepsilon,\delta>0$, and find $c>0$ and $a\in\left(  0,c\right)  $ such that
$\delta<am_{a}$ and $m_{a}<\varepsilon.$ Looking at the inequality concerning
$\delta,$ it is important to maximize the function
\[
\left(  0,c\right)  \ni a\mapsto a\left(  \dfrac{c}{c+1}-\dfrac{a}%
{a+1}\right),
\]
and an elementary computation shows that $\overline{a}=\sqrt{1+c}-1$ is the
maximum point. For this choice of $a,$ we have
\[
m_{\overline{a}}=\min\left\{  c,\frac{1}{4}\left(  \dfrac{1}{\sqrt{c+1}%
}-\dfrac{1}{c+1}\right)  ,c\left(  \frac{1}{\sqrt{c+1}}+\dfrac{1}{2\left(
c+1\right)  }\right)  \right\}  .
\]
Now, it is not very involved to show that, actually,
\[
m_{\overline{a}}=\frac{1}{4}\left(  \dfrac{1}{\sqrt{c+1}}-\dfrac{1}%
{c+1}\right)  .
\]
So, we have to find $c>0$ such that the following system of inequalities
holds:%
\begin{align}
&  \delta<\frac{1}{4}\left(  \sqrt{1+c}-1\right)  \left(  \dfrac{1}{\sqrt
{c+1}}-\dfrac{1}{c+1}\right),  \label{sist_ineq}\\
&  \frac{1}{4}\left(  \dfrac{1}{\sqrt{c+1}}-\dfrac{1}{c+1}\right)
<\varepsilon.\nonumber
\end{align}
This is a slight variation of a system of inequalities discussed in
\cite{D2020} and \cite{C2022}. So, according to \cite[Lemma 3.2]{C2022} (with
some obvious changes), if $\varepsilon\in\left(  0,\frac{1}{16}\right)  $, this
system has a solution provided $\sqrt{\delta}-\delta<\varepsilon$ and
$\delta<\varepsilon$, and for $\varepsilon>\frac{1}{16}$, the system has a
solution if and only if $\delta\in\left(  0,\frac{1}{4}\right)  .$ Of course,
having in view the fact that the inequalities are strict and all the
underlying functions are continuous, once the system has a solution, it
actually has as solutions all the points in an interval.

Notice that the function $\varphi:[0,\infty)\rightarrow\left[  0,\frac{1}%
{4}\right)  $ given by%
\begin{equation}
\varphi\left(  x\right)  :=\frac{1}{4}\left(  \sqrt{1+x}-1\right)  \left(
\dfrac{1}{\sqrt{x+1}}-\dfrac{1}{x+1}\right)  \label{fct_phi}%
\end{equation}
is an increasing continuous bijection.

Denote $\psi:[0,\infty)\rightarrow\mathbb{R}$,%
\[
\psi\left(  x\right)  :=\frac{1}{4}\left(  \dfrac{1}{\sqrt{x+1}}-\dfrac
{1}{x+1}\right)  .
\]

Let us remind that {a} set-valued {map} $F$ is said to have the Aubin property
around $(\overline{x},\overline{y})$ with constant $\ell_{F}>0$ if there exist
two neighborhoods $U$ of $\overline{x}$ and $V\ $of $\overline{y}$ such that,
for every $x,u\in U,$%
\[
F(x)\cap V\subset F(u)+\ell_{F}\left\Vert x-u\right\Vert \overline{B}_{Y}.
\]
It became classical that the latter property is equivalent to the openness at
linear rate for $F^{-1}$ around $(\overline{y},\overline{x})$. For more
details, see, for instance, \cite{AubFra}, \cite{DR2009}. Moreover, it also
became classical the following characterization of Aubin property by the
Mordukhovich coderivative (see, e.g., \cite{M2018}).

\begin{lm}
[Mordukhovich Criterion]\label{MC}Let $X,Y$ be finite dimensional spaces,
$G:X\rightrightarrows Y$ be a closed-graph set-valued mapping and $(\overline
{x},\overline{y})\in\operatorname*{Gr}G.$ Then $G$ has Aubin property around
$(\overline{x},\overline{y})$ iff $D^{\ast}G\left(  \overline{x},\overline
{y}\right)  \left(  0\right)  =\left\{  0\right\}  .$
\end{lm}

The next result contains the necessary conditions we envisaged.

\begin{thm}
\label{nec_cond_set-valued_K}Let $X,Y$ be finite dimensional spaces,
$\Omega\subseteq X$ be a closed set, $F:X\rightrightarrows Y$ be a
closed-graph set-valued mapping, $K:X\rightrightarrows Y$ be a closed-graph vds,
and $\left(  \overline{x},\overline{y}\right)  \in X\times Y$ such that
$\overline{x}\in\Omega,(\overline{x},\overline{y})\in\operatorname*{Gr}%
F,\left(  \overline{x},0\right)  \in\operatorname*{Gr}K.$ Moreover, suppose
that the following conditions are satisfied:

(a) (\textsc{Aubin Conditions}) $F$ and $K$ have Aubin property around
$(\overline{x},\overline{y})$ and $\left(  \overline{x},0\right)  ,$ respectively.

(b) (\textsc{Constants Conditions}) There exists $\varepsilon\in\left(
0,\infty\right]  $ and $\delta\in\left(  0,\infty\right)  $ such that one of
the following set of conditions hold:

$\quad$(b.i) $\varepsilon>\frac{1}{16}$ and $\delta\in\left(  0,\frac{1}%
{4}\right)  ;$

$\quad$(b.ii) $\delta<\varepsilon$ and $\sqrt{\delta}-\delta<\varepsilon.$

(c) (\textsc{Interiority Condition}) For $\varepsilon$ from above, denote
$P:=%
{\displaystyle\bigcap_{x\in B\left(  \overline{x},\varepsilon\right)  }}
K\left(  x\right)  $ and suppose that $k\in S_{Y}\cap\operatorname*{int}P.$

(d) (\textsc{Approximate Solution Condition}) $\left(  \overline{x}%
,\overline{y}\right)  $ is an $\left(  \varepsilon,\delta,{{k}}\right)
${-nondominated solution} wrt $K$ for $F$ on $\Omega$ .

\noindent Let $\varphi$ be the function given by (\ref{fct_phi}). Then, there
exist $y^{\ast}\in P^{+}$ with $\left\langle y^{\ast},k\right\rangle =1,$ and
$(x_{1},y_{1})\in\operatorname*{Gr}F\cap\lbrack\overline{B}(\overline
{x},\varphi^{-1}\left(  \delta\right)  )\times\overline{B}(\overline
{y},\varphi^{-1}\left(  \delta\right)  )],$ $(x_{2},y_{2})\in
\operatorname*{Gr}K\cap\lbrack\overline{B}(\overline{x},\varphi^{-1}\left(
\delta\right)  )\times\overline{B}(0,\varphi^{-1}\left(  \delta\right)  )],$
$x_{3}\in\Omega\cap\overline{B}\left(  \overline{x},\varphi^{-1}\left(
\delta\right)  \right)  ,$ $z^{\ast}\in2\varphi^{-1}\left(  \delta\right)
\overline{B}_{Y^{\ast}},$ such that%
\[
0\in D^{\ast}F(x_{1},y_{1})(y^{\ast}-z^{\ast})+D^{\ast}K(x_{2},y_{2})(y^{\ast
})+N\left(  \Omega,x_{3}\right)  +\varphi^{-1}\left(  \delta\right)
\overline{B}_{X^{\ast}}.
\]

\end{thm}

\noindent\textbf{Proof. }Using one of the sets of assumptions (b.i) and
(b.ii), we get by the discussion above that the system of inequalities
(\ref{sist_ineq}) has an interval of solutions. Because $\varphi$ is a
increasing continuous bijection, we can take this interval of the form
$\left(  \varphi^{-1}\left(  \delta\right)  ,\varphi^{-1}\left(
\delta\right)  +\mu\right)  ,$ with $\mu>0.$ Moreover, since $\delta
<\varphi\left(  c\right)  \ $is equivalent to $\varphi^{-1}\left(
\delta\right)  <c,$ we obtain that $B(0,\varphi^{-1}\left(  \delta\right)
)\subseteq B\left(  0,c\right)  .$

Fix $d\in\left(  \varphi^{-1}\left(  \delta\right)  ,\varphi^{-1}\left(
\delta\right)  +\mu\right)  $ (hence $\delta<\varphi\left(  d\right)  $) and
take $\rho:=\varphi\left(  d\right)  ,$ hence $d=\varphi^{-1}\left(
\rho\right)  .$ Observe that, for any $\delta^{\prime}\in\left(  \delta
,\rho\right)  $ sufficiently close to $\delta,$ we have by continuity that%
\[
\psi\left(  \varphi^{-1}\left(  \delta^{\prime}\right)  \right)  =\frac{1}%
{4}\left(  \dfrac{1}{\sqrt{\varphi^{-1}\left(  \delta^{\prime}\right)  +1}%
}-\dfrac{1}{\varphi^{-1}\left(  \delta^{\prime}\right)  +1}\right)
<\varepsilon.
\]
Hence, for any $\delta^{\prime}\in\left(  \delta,\rho\right)  $ sufficiently
close to $\delta,$ the system of inequalitites%
\begin{align}
&  \delta<\frac{1}{4}\left(  \sqrt{1+\varphi^{-1}\left(  \delta^{\prime
}\right)  }-1\right)  \left(  \dfrac{1}{\sqrt{\varphi^{-1}\left(
\delta^{\prime}\right)  +1}}-\dfrac{1}{\varphi^{-1}\left(  \delta^{\prime
}\right)  +1}\right)  =\delta^{\prime}\label{sist_ineq2}\\
&  \frac{1}{4}\left(  \dfrac{1}{\sqrt{\varphi^{-1}\left(  \delta^{\prime
}\right)  +1}}-\dfrac{1}{\varphi^{-1}\left(  \delta^{\prime}\right)
+1}\right)  <\varepsilon\nonumber
\end{align}
is satisfied.

Take $\delta_{n}\downarrow\delta$ and consider, for any $n,$ $\delta
_{n}^{\prime}\in\left(  \delta,\delta_{n}\right)  .$ It follows that%
\begin{align*}
&  \delta<\delta_{n}^{\prime}<\delta_{n}=\frac{1}{4}\left(  \sqrt
{1+\varphi^{-1}\left(  \delta_{n}\right)  }-1\right)  \left(  \dfrac{1}%
{\sqrt{\varphi^{-1}\left(  \delta_{n}\right)  +1}}-\dfrac{1}{\varphi
^{-1}\left(  \delta_{n}\right)  +1}\right)  ,\\
&  \frac{1}{4}\left(  \dfrac{1}{\sqrt{\varphi^{-1}\left(  \delta_{n}\right)
+1}}-\dfrac{1}{\varphi^{-1}\left(  \delta_{n}\right)  +1}\right)
<\varepsilon.
\end{align*}
We have, for%
\begin{align*}
a_{n}  &  :=\sqrt{1+\varphi^{-1}\left(  \delta_{n}\right)  }-1,\\
\rho_{n}  &  :=\psi\left(  \varphi^{-1}\left(  \delta_{n}\right)  \right)  ,
\end{align*}
that $a_{n}\rho_{n}=\delta_{n}>\delta_{n}^{\prime}$ and $\rho_{n}%
<\varepsilon.$

By Proposition \ref{prop_incom}, we know that%
\[
B\left(  \overline{y},\delta_{n}^{\prime}\right)  \cap\left[  \overline
{y}-\operatorname*{cone}\left\{  k\right\}  \right]  \nsubseteq\left(
F+K+\Delta_{\Omega}\right)  \left(  B\left(  \overline{x},\varepsilon\right)
\right)  ,
\]
which means that%
\[
B\left(  \overline{y},a_{n}\rho_{n}\right)  \cap\left[  \overline
{y}-\operatorname*{cone}\left\{  k\right\}  \right]  \nsubseteq\left(
F+K+\Delta_{\Omega}\right)  \left(  B\left(  \overline{x},\rho_{n}\right)
\right)  ,
\]
which shows that the conclusion of Theorem \ref{dir_op_FKD} cannot be
satisfied for any $n.$

Remark that, using (\textsc{Aubin Conditions}) and Lemma \ref{MC}, we obtain
that $D^{\ast}F\left(  \overline{x},\overline{y}\right)  \left(  0\right)
=\left\{  0\right\}  $ and $D^{\ast}K\left(  \overline{x},0\right)  \left(
0\right)  =\left\{  0\right\}  ,$ hence the (\textsc{Transversality
Conditions}) from Theorem \ref{dir_op_FKD} are satisfied.

Hence, all the assumptions of the Theorem, except (\textsc{Injectivity
Condition}), are satisfied, hence we deduce that for any $n,$ there exist
$y_{n}^{\ast}\in P^{+}$ with $\left\langle y_{n}^{\ast},k\right\rangle =1,$
and $(x_{1n},y_{1n})\in\operatorname*{Gr}F\cap\lbrack B(\overline{x}%
,\varphi^{-1}\left(  \delta_{n}\right)  )\times B(\overline{y},\varphi
^{-1}\left(  \delta_{n}\right)  )],$ $(x_{2n},y_{2n})\in\operatorname*{Gr}%
K\cap\lbrack B(\overline{x},\varphi^{-1}\left(  \delta_{n}\right)  )\times
B(0,\varphi^{-1}\left(  \delta_{n}\right)  )],$ $x_{3n}\in\Omega\cap B\left(
\overline{x},\varphi^{-1}\left(  \delta_{n}\right)  \right)  ,$ $z_{n}^{\ast
}\in2\varphi^{-1}\left(  \delta_{n}\right)  B_{Y^{\ast}},$
$x_{1n}^{\ast}\in\widehat{D}^{\ast}F(x_{1n},y_{1n})(y_{n}^{\ast}-z_{n}^{\ast
}),$ $x_{2n}^{\ast}\in\widehat{D}^{\ast}K(x_{2n},y_{2n})(y_{n}^{\ast}),$
$x_{3n}^{\ast}\in\widehat{N}\left(  \Omega,x_{3n}\right)  $ such that%
\begin{equation}
\varphi^{-1}\left(  \delta_{n}\right)  \geq\left\Vert x_{1n}^{\ast}%
+x_{2n}^{\ast}+x_{3n}^{\ast}\right\Vert .\label{ineq}%
\end{equation}
Using $k\in S_{Y}\cap\operatorname*{int}P$ and the fact (see \cite[Corollary
2.2.35]{GRTZ2023}) that the set%
\[
\left\{  y^{\ast}\in P^{+}\mid\left\langle y^{\ast},k\right\rangle =1\right\}
\]
is a (weak$^{\ast}$) base for $P^{+},$ we deduce that $\left(  y_{n}^{\ast
}\right)  $ is bounded. Moreover, since $F$ and $K$ have Aubin property, we
deduce that $\left(  x_{1n}^{\ast}\right)  $ and $\left(  x_{2n}^{\ast
}\right)  $ are bounded, hence from (\ref{ineq}), $\left(  x_{3n}^{\ast
}\right)  $ is also bounded. Since the sequences $(x_{1n},y_{1n}%
),(x_{2n},y_{2n})$ and $(x_{3n})$ are also bounded, and all the involved
spaces are finite dimensional, we may suppose, without relabeling, that all
of them converge to corresponding points. Moreover, using also that the
$\varphi^{-1}$ is continuous, we obtain that $y_{n}^{\ast}\rightarrow y^{\ast
}\in P^{+}$ (since $P^{+}$ is closed) and $\left\langle y^{\ast}%
,k\right\rangle =1,$ $z_{n}^{\ast}\rightarrow z^{\ast}\in2\varphi^{-1}\left(
\delta\right)  \overline{B}_{Y^{\ast}},$ $(x_{1n},y_{1n}%
)\rightarrow\left(  x_{1},y_{1}\right)  \in\operatorname*{Gr}F\cap
\lbrack\overline{B}(\overline{x},\varphi^{-1}\left(  \delta\right)
)\times\overline{B}(\overline{y},\varphi^{-1}\left(  \delta\right)  )]$
($\operatorname*{Gr}F$ is closed, $\varphi^{-1}$ is continuous),
$(x_{2n},y_{2n})\rightarrow(x_{2},y_{2})\in\operatorname*{Gr}K\cap
\lbrack\overline{B}(\overline{x},\varphi^{-1}\left(  \delta\right)
)\times\overline{B}(0,\varphi^{-1}\left(  \delta\right)  )]$
($\operatorname*{Gr}K$ is closed, $\varphi^{-1}$ is continuous),
$x_{3n}\rightarrow x_{3}\in\Omega\cap\overline{B}(\overline{x},\varphi
^{-1}\left(  \delta\right)  )$ ($\Omega$ is closed, $\varphi^{-1}$ is
continuous), $x_{1n}^{\ast}\rightarrow x_{1}^{\ast}\in D^{\ast}F(x_{1}%
,y_{1})(y^{\ast}-z^{\ast}),$ $x_{2n}^{\ast}\rightarrow x_{2}^{\ast}\in
D^{\ast}K(x_{2},y_{2})(y^{\ast}),$ $x_{3n}^{\ast}\rightarrow x_{3}^{\ast}\in
N\left(  \Omega,x_{3}\right)  .$ Moreover, we get from (\ref{ineq}) that
$\varphi^{-1}\left(  \delta\right)  \geq\left\Vert x_{1}^{\ast}+x_{2}^{\ast
}+x_{3}^{\ast}\right\Vert .$ Therefore,%
\[
0\in D^{\ast}F(x_{1},y_{1})(y^{\ast}-z^{\ast})+D^{\ast}K(x_{2}%
,y_{2})(y^{\ast})+N\left(  \Omega,x_{3}\right)  +\varphi^{-1}\left(
\delta\right)  \overline{B}_{X^{\ast}},
\]
which shows that the conclusion is satisfied.$\hfill\square$

\subsection{Necessary conditions for {  approximately nondominated} points wrt $Q$}

Next, we want to obtain a similar result for the case of domination wrt
$Q:Y\rightrightarrows Y.$ By Remark \ref{reduction}, since $\left(
\overline{x},\overline{y}\right)  $ is an $\left(  \varepsilon,\delta,{{k}%
}\right)  $-nondominated solution wrt $Q$ for $F$ on $\Omega$, it follows
that $\left(  \left(  \overline{x},\overline{y}\right)  ,\overline{y}\right)
\in\operatorname*{Gr}\widetilde{F}$ is an $\left(  \left(  \varepsilon
,+\infty\right)  ,\delta,{{k}}\right)  $-nondominated solution wrt
$\widetilde{K}$ for $\widetilde{F}$ on $\Omega\times Y$, where $\widetilde{F}$
and $\widetilde{K}$ are given by (\ref{F_tilt}) and (\ref{K_tilt}), respectively.

The next lemma shows the links between the coderivatives of the involved
set-valued  mappings (see \cite{DST2024}).

\begin{lm}
Let {$X,Y$ be Asplund spaces}, $F:X\rightrightarrows Y,Q:Y\rightrightarrows Y$
be set-valued mappings, and suppose that $\widetilde{F},\widetilde{K}:X\times
Y\rightrightarrows Y$ are given by (\ref{F_tilt}), (\ref{K_tilt}). We have:%
\begin{align}
\left(  x^{\ast},y^{\ast}\right)   &  \in\widehat{D}^{\ast}\widetilde
{F}(x,y,y)(z^{\ast})\Leftrightarrow x^{\ast}\in\widehat{D}^{\ast
}F(x,y)(z^{\ast}-y^{\ast}),\label{Frcoder_F_tilt}\\
\left(  x^{\ast},y^{\ast}\right)   &  \in\widehat{D}^{\ast}\widetilde
{K}(x,y,z)(z^{\ast})\Leftrightarrow x^{\ast}=0,y^{\ast}\in\widehat{D}^{\ast
}Q(y,z)(z^{\ast}). \label{Frcoder_K_tilt}%
\end{align}
Moreover, if $F$ and $Q$ are closed-graph, then $\widetilde{F},\widetilde{K}$
are also closed-graph, and%
\begin{align}
\left(  x^{\ast},y^{\ast}\right)   &  \in D^{\ast}\widetilde{F}(x,y,y)(z^{\ast
})\Leftrightarrow x^{\ast}\in D^{\ast}F(x,y)(z^{\ast}-y^{\ast}%
),\label{Mcoder_F_tilt}\\
\left(  x^{\ast},y^{\ast}\right)   &  \in D^{\ast}\widetilde{K}(x,y,z)(z^{\ast
})\Leftrightarrow x^{\ast}=0,y^{\ast}\in D^{\ast}Q(y,z)(z^{\ast}).
\label{Mcoder_K_tilt}%
\end{align}

\end{lm}

Now, it would be natural to try to apply Theorem \ref{nec_cond_set-valued_K}.
Unfortunately, one can see that $F$ has Aubin property around $(\overline
{x},\overline{y})$ iff $D^{\ast}F\left(  \overline{x},\overline{y}\right)
\left(  0\right)  =\left\{  0\right\}  $ and, using (\ref{Mcoder_F_tilt}),
this only implies that $\left(  0,y^{\ast}\right)  \in D^{\ast}\widetilde
{F}(x,y,y)(y^{\ast})$ for $y^{\ast}\in Y^{\ast},$ hence $\widetilde{F}$
doesn't have Aubin property at $\left(  \overline{x},\overline{y},\overline
{y}\right)  .$ Moreover, given the form of the set-valued $\widetilde{F},$ it
would seem very unnatural to impose its Aubin property.

Nevertheless, by Proposition \ref{prop_incom}, we know that, in this case, if
\[
q\in\left(
{\displaystyle\bigcap_{\left(  x,y\right)  \in\left[  B\left(  \overline
{x},\varepsilon\right)  \cap\Omega\right]  \times B\left(  \overline
{y},\varepsilon\right)  }}
\widetilde{K}\left(  x,y\right)  \right)  \setminus\left\{  0\right\}
=\left(
{\displaystyle\bigcap_{y\in B\left(  \overline{y},\varepsilon\right)  }}
Q\left(  y\right)  \right)  \setminus\left\{  0\right\}  ,
\]
and $\left(  \left(  \overline{x},\overline{y}\right)  ,\overline{y}\right)
\in\operatorname*{Gr}\widetilde{F}$ is an $\left(  \left(  \varepsilon
,+\infty\right)  ,\delta,{{k}}\right)  ${-nondominated solution} wrt
$\widetilde{K}$ for $\widetilde{F}$ on $\Omega\times Y,$ then one cannot have
$\delta^{\prime}>\delta$ such that%
\[
B\left(  \overline{y},\delta^{\prime}\right)  \cap\left[  \overline
{y}-\operatorname*{cone}\left\{  q\right\}  \right]  \subseteq\left(
\widetilde{F}+\widetilde{K}+\Delta_{\Omega\times Y}\right)  \left(  B\left(
\left(  \overline{x},\overline{y}\right)  ,\varepsilon\right)  \right)  .
\]

So, we are forced to go back to Theorem \ref{suf_cond_dir_op_sum}. First, we
want to check when the sets corresponding to (\ref{D1p}) for the mappings
$\widetilde{F},\widetilde{K}$ and $\Delta_{\Omega\times Y}$ are allied at
$\left(  \overline{x},\overline{y},\overline{y},0,0\right)  $. These sets are%
\begin{align}
E_{1}  &  =\{(x,y,z,w,v)\in X\times Y^{4}\mid z\in\widetilde{F}%
(x,y)\}=\{(x,y,y,w,v)\in X\times Y^{4}\mid y\in F(x)\},\nonumber\\
E_{2}  &  =\{(x,y,z,w,v)\in X\times Y^{4}\mid w\in\widetilde{K}%
(x,y)\}=\{(x,y,z,w,v)\in X\times Y^{4}\mid w\in Q(y)\},\label{E123}\\
E_{3}  &  =\{(x,y,z,w,v)\in X\times Y^{4}\mid v\in\Delta_{\Omega\times
Y}(x,y)\}=\{(x,y,z,w,0)\in X\times Y^{4}\mid x\in\Omega\}.\nonumber
\end{align}
One can easily check that the next assertion holds (see \cite{DST2024}).

\begin{lm}
\label{alliedness_E_imply_A}The sets $E_{1},E_{2},E_{3}$ are allied at
$\left(  \overline{x},\overline{y},\overline{y},0,0\right)  $ {if and only if}
the sets%
\begin{align}
A_{1}  &  :=\{(x,y,z)\in X\times Y^{2}\mid y\in F(x)\},\nonumber\\
A_{2}  &  :=\{(x,y,z)\in X\times Y^{2}\mid z\in Q(y)\},\label{A13}\\
A_{3}  &  :=\{(x,y,z)\in X\times Y^{2}\mid x\in\Omega\}\nonumber
\end{align}
are allied at $\left(  \overline{x},\overline{y},0\right)  .$
\end{lm}

Moreover, similar to the previous case, some transversality conditions imposed
on $F,Q$ and $\Omega$ imply the alliedness of these sets (see \cite{DST2024}).

\begin{lm}
\label{alliedness2}Let $X,Y$ be Asplund spaces, $\Omega\subseteq X$ be a
closed set, $F:X\rightrightarrows Y,Q:Y\rightrightarrows Y$ be closed-graph
set-valued mappings, and $\left(  \overline{x},\overline{y},0\right)  \in
X\times Y\times Y$ such that $\overline{x}\in\Omega,(\overline{x},\overline
{y})\in\operatorname*{Gr}F,\left(  \overline{y},0\right)  \in
\operatorname*{Gr}Q.$ Moreover, assume that:%
\begin{align*}
&  D^{\ast}Q(\overline{y},0)(0)=\left\{  0\right\}  ,\\
&  D^{\ast}F(\overline{x},\overline{y})(0)\cap\left(  -N\left(  \Omega
,\overline{x}\right)  \right)  =\left\{  0\right\}  .
\end{align*}
Then the sets $A_{1},A_{2},A_{3}$ defined by (\ref{A13}) are allied at
$\left(  \overline{x},\overline{y},0\right)  .$
\end{lm}

Remark also that the sum set-valued  mapping $H:=\widetilde{F}+\widetilde{K}%
+\Delta_{\Omega\times Y}$ given by (\ref{H_sum})\ becomes in this particular
case%
\begin{equation}
H\left(  x,y\right)  =\left\{
\begin{array}
[c]{ll}%
\left\{  y\right\}  +Q\left(  y\right)  & \text{if }y\in F\left(  x\right)
,x\in\Omega\\
\emptyset & \text{otherwise.}%
\end{array}
\right.  \label{Hxy}%
\end{equation}

We are now in position to give a directional openness result of the type of
Theorem \ref{dir_op_FKD}.

\begin{thm}
\label{dir_op_FQD}Let $X,Y$ be finite dimensional spaces, $\Omega\subseteq X$
be a closed set, $F:X\rightrightarrows Y$ be a closed-graph set-valued mapping,
$Q:Y\rightrightarrows Y$ be a closed-graph vds, and $\left(  \overline
{x},\overline{y}\right)  \in X\times Y$ such that $\overline{x}\in
\Omega,(\overline{x},\overline{y})\in\operatorname*{Gr}F,\left(  \overline
{y},0\right)  \in\operatorname*{Gr}Q.$ Denote $P:=%
{\displaystyle\bigcap_{y\in B\left(  \overline{y},\eta\right)  }}
Q\left(  y\right)  $ for some $\eta>0$ and suppose $q\in P\cap S_{Y}.$
Moreover, suppose that the following conditions are satisfied:

(a) (\textsc{Transversality Conditions})
\begin{align*}
&  D^{\ast}Q(\overline{y},0)(0)=\left\{  0\right\}  ,\\
&  D^{\ast}F(\overline{x},\overline{y})(0)\cap\left(  -N\left(  \Omega
,\overline{x}\right)  \right)  =\left\{  0\right\}  .
\end{align*}

(b) (\textsc{Injectivity Condition}) there exist $r\in\left(  0,\eta\right)
,c>0$ such that for every $y^{\ast}\in P^{+}$ with $\left\langle y^{\ast
},q\right\rangle =1,$ every $\left(  {x_{1},y_{1}}\right)  {\in}%
\operatorname*{Gr}F\cap\lbrack B(\overline{x},r)\times B(\overline{y},r)],$
${(}y_{2}{,z_{2})\in}\operatorname*{Gr}Q\cap\lbrack B(\overline{y},r)\times
B(0,r)],$ ${x_{3}\in\Omega}\cap B(\overline{x},r),$ every $y_{1}^{\ast}\in
Y^{\ast},z^{\ast}\in2cB_{Y^{\ast}},$ and every $x_{1}^{\ast}\in\widehat
{D}^{\ast}F({x_{1},y_{1}})(y^{\ast}-y_{1}^{\ast}-z^{\ast}),$ $y_{2}^{\ast}%
\in\widehat{D}^{\ast}Q({y_{2},z_{2}})(y^{\ast}),$ $x_{3}^{\ast}\in\widehat
{N}\left(  \Omega,x_{3}\right)  ,$ one has%
\begin{equation}
c\leq\Vert\left(  x_{1}^{\ast}+x_{3}^{\ast},y_{1}^{\ast}+y_{2}^{\ast}\right)
\Vert.\label{ineqxy2}%
\end{equation}

\noindent Take arbitrary $a\in(0,c),$ and $\theta$ such that (\ref{epsilon})
holds. Then for any $\rho\in\left(  0,\theta\right)  ,$%
\[
B(\overline{y},a\rho)\cap\left[  \overline{y}-\operatorname*{cone}\left\{
q\right\}  \right]  \subseteq H\left(  B(\left(  \overline{x},\overline
{y}\right)  ,\rho)\right)  .
\]

\end{thm}

\textbf{\noindent Proof.} By Lemmata \ref{alliedness_E_imply_A} and
\ref{alliedness2}, we see that the (\textsc{Transversality Conditions}) imply
the alliedness of $E_{1},E_{2},E_{3}$ defined by (\ref{E123}) at $\left(
\overline{x},\overline{y},\overline{y},0,0\right)  .$ Moreover, remark that,
besides (\ref{Frcoder_F_tilt}) and (\ref{Frcoder_K_tilt}), we also have, for
any $x\in\Omega,$%
\begin{equation}
\left(  x^{\ast},y^{\ast}\right)  \in\widehat{D}^{\ast}\Delta_{\Omega\times
Y}(x,y,0)(v^{\ast})\Leftrightarrow x^{\ast}\in\widehat{N}\left(
\Omega,x\right)  ,y^{\ast}=0,v^{\ast}\in Y^{\ast}.\label{Frcoder_D}%
\end{equation}
We want to check how the final assumption of Theorem \ref{suf_cond_dir_op_sum}
looks in our case.

We should find $c>0$ and $r>0$ such that for every ${(}\left(  {x_{1},y_{1}%
}\right)  {,y_{1})\in}\operatorname*{Gr}\widetilde{F}\cap\lbrack
B(\overline{x},r)\times B(\overline{y},r)\times B(\overline{y},r)],$
${(}\left(  {x_{2},y_{2}}\right)  {,z_{2})\in}\operatorname*{Gr}\widetilde
{K}\cap\lbrack B(\overline{x},r)\times B(\overline{y},r)\times B(0,r)],$
${(}\left(  {x_{3},y_{3}}\right)  {,0)\in}\operatorname*{Gr}\Delta
_{\Omega\times Y}\cap\lbrack B(\overline{x},r)\times B(\overline{y},r)\times
B(0,r)]$ every $y^{\ast}\in Y^{\ast},$ with $\left\langle y^{\ast
},q\right\rangle =1,z_{i}^{\ast}\in2cB_{Y^{\ast}},i=\overline{{1,2}}$ and
every $\left(  x_{1}^{\ast},y_{1}^{\ast}\right)  \in\widehat{D}^{\ast
}\widetilde{F}(\left(  {x_{1},y_{1}}\right)  {,y_{1}})(y^{\ast}-z_{1}^{\ast
}),$ $\left(  x_{2}^{\ast},y_{2}^{\ast}\right)  \in\widehat{D}^{\ast
}\widetilde{K}(\left(  {x_{2},y_{2}}\right)  {,z_{2}})(y^{\ast}),$ $\left(
x_{3}^{\ast},y_{3}^{\ast}\right)  \in\widehat{D}^{\ast}\Delta_{\Omega\times
Y}(\left(  {x_{3},y_{3}}\right)  {,0})(y^{\ast}-z_{2}^{\ast}),$ one has%
\begin{equation}
c\leq\Vert\left(  x_{1}^{\ast},y_{1}^{\ast}\right)  +\left(  x_{2}^{\ast
},y_{2}^{\ast}\right)  +\left(  x_{3}^{\ast},y_{3}^{\ast}\right)  \Vert.
\label{ineqxy}%
\end{equation}

This reduces, by the use of (\ref{Frcoder_F_tilt}), (\ref{Frcoder_K_tilt}) and
(\ref{Frcoder_D}), exactly to the (\textsc{Injectivity Condition}) above, if
we prove that necessarily the element $y^{\ast}$ involved in the condition is
from $P^{+}.$ But this follows in a similar way as in the proof of Theorem
\ref{dir_op_FKD}.$\hfill\square$

\bigskip

Now, we are in position to formulate the necessary conditions for
approximate minima wrt vds $Q$.

\begin{thm}
\label{nec_cond_set-valued_Q}Let $X,Y$ be finite dimensional spaces,
$\Omega\subseteq X$ be a closed set, $F:X\rightrightarrows Y$ be a
closed-graph set-valued mapping, $Q:Y\rightrightarrows Y$ be a closed-graph vds
and $\left(  \overline{x},\overline{y}\right)  \in X\times Y$ such that
$\overline{x}\in\Omega,(\overline{x},\overline{y})\in\operatorname*{Gr}%
F,\left(  \overline{y},0\right)  \in\operatorname*{Gr}Q.$ Moreover, suppose
that the following conditions are satisfied:

(a) (\textsc{Aubin Conditions}) $F$ and $Q$ have Aubin property around
$(\overline{x},\overline{y})$ and $\left(  \overline{y},0\right)  ,$ respectively.

(b) (\textsc{Constants Conditions}) There exists $\varepsilon\in\left(
0,\infty\right]  $ and $\delta\in\left(  0,\infty\right)  $ such that one of
the following set of conditions hold:

$\quad$(b.i) $\varepsilon>\frac{1}{16}$ and $\delta\in\left(  0,\frac{1}%
{4}\right)  ;$

$\quad$(b.ii) $\delta<\varepsilon$ and $\sqrt{\delta}-\delta<\varepsilon.$

(c) (\textsc{Interiority Conditions}) For $\varepsilon$ from above, denote
$P:=%
{\displaystyle\bigcap_{y\in B\left(  \overline{y},\varepsilon\right)  }}
Q\left(  y\right)  $ and suppose that $q\in S_{Y}\cap\operatorname*{int}P.$

(d) (\textsc{Approximate Solution Condition}) $\left(  \overline{x}%
,\overline{y}\right)  $ is an $\left(  \varepsilon,\delta,{{k}}\right)
${-nondominated solution} wrt $Q$ for $F$ on $\Omega$ .

\noindent Let $\varphi$ be the function given by (\ref{fct_phi}). Then, there
exist $y^{\ast}\in P^{+}$ with $\left\langle y^{\ast},q\right\rangle =1,$ and
$(x_{1},y_{1})\in\operatorname*{Gr}F\cap\lbrack\overline{B}(\overline
{x},\varphi^{-1}\left(  \delta\right)  )\times\overline{B}(\overline
{y},\varphi^{-1}\left(  \delta\right)  )],$ $(y_{2},z_{2})\in
\operatorname*{Gr}Q\cap\lbrack\overline{B}(\overline{y},\varphi^{-1}\left(
\delta\right)  )\times\overline{B}(0,\varphi^{-1}\left(  \delta\right)  )],$
$x_{3}\in\Omega\cap\overline{B}\left(  \overline{x},\varphi^{-1}\left(
\delta\right)  \right)  ,$ $z^{\ast}\in2\varphi^{-1}\left(  \delta\right)
\overline{B}_{Y^{\ast}},$ and $y_{1}^{\ast}\in\widehat{D}^{\ast}Q({y_{2}%
,z_{2}})(y^{\ast})+\varphi^{-1}\left(  \delta\right)  \overline{B}_{Y^{\ast}}$
such that%
\[
0\in D^{\ast}F(x_{1},y_{1})(y^{\ast}+y_{1}^{\ast}+z^{\ast})+N\left(
\Omega,x_{3}\right)  +\varphi^{-1}\left(  \delta\right)  \overline{B}%
_{X^{\ast}}.
\]

\end{thm}

\noindent\textbf{Proof.} Remark that the Aubin conditions imposed on $F$ and
$Q$ imply that the (\textsc{Transversality Conditions}) from Theorem
\ref{dir_op_FQD} hold. Making the same steps as in the proof of Theorem
\ref{nec_cond_set-valued_K}, we arrive that the (\textsc{Injectivity
Condition}) from Theorem \ref{dir_op_FQD} cannot be satisfied. We deduce that,
taking $\delta_{n}\downarrow\delta,$ for any $n,$ there exist $y_{n}^{\ast}\in
P^{+}$ with $\left\langle y_{n}^{\ast},q\right\rangle =1,$ and $(x_{1n}%
,y_{1n})\in\operatorname*{Gr}F\cap\lbrack B(\overline{x},\varphi^{-1}\left(
\delta_{n}\right)  )\times B(\overline{y},\varphi^{-1}\left(  \delta
_{n}\right)  )],$ $(y_{2n},z_{2n})\in\operatorname*{Gr}Q\cap\lbrack
B(\overline{y},\varphi^{-1}\left(  \delta_{n}\right)  )\times B(0,\varphi
^{-1}\left(  \delta_{n}\right)  )],$ $x_{3n}\in\Omega\cap B\left(
\overline{x},\varphi^{-1}\left(  \delta_{n}\right)  \right)  ,$ $y_{1n}^{\ast
}\in Y^{\ast}, z_{n}^{\ast}\in2\varphi^{-1}\left(  \delta_{n}\right)
B_{Y^{\ast}},$ and $x_{1n}^{\ast}\in\widehat{D}^{\ast}F({x_{1n},y_{1n}}%
)(y_{n}^{\ast}-y_{1n}^{\ast}-z_{n}^{\ast}),$ $y_{2n}^{\ast}\in\widehat
{D}^{\ast}Q({y_{2n},z_{2n}})(y_{n}^{\ast}),$ $x_{3n}^{\ast}\in\widehat
{N}\left(  \Omega,x_{3}\right)  ,$ such that%
\begin{equation}
\varphi^{-1}\left(  \delta_{n}\right)  \geq\Vert\left(  x_{1n}^{\ast}%
+x_{3n}^{\ast},y_{1n}^{\ast}+y_{2n}^{\ast}\right)  \Vert.\label{ineqxyn}%
\end{equation}
Using $q\in S_{Y}\cap\operatorname*{int}P,$ we deduce again that $\left(
y_{n}^{\ast}\right)  $ is bounded. Since $Q$ has Aubin property, it follows
that $\left(  y_{2n}^{\ast}\right)  $ is bounded, hence by (\ref{ineqxyn}), it
follows that $\left(  y_{1n}^{\ast}\right)  $ is bounded. This shows that the
sequence $(y_{n}^{\ast}-y_{1n}^{\ast}-z_{1n}^{\ast})$ is bounded. But $F$ has
Aubin property, so also $\left(  x_{1n}^{\ast}\right)  $ must be bounded, so
finally, using again (\ref{ineqxyn}), it follows that $\left(  x_{3n}^{\ast
}\right)  $ is bounded. The proof ends similarly to the proof of Theorem
\ref{nec_cond_set-valued_K}.$\hfill\square$

\section{Conclusions}\label{s-concl}
{  In our paper, we derived new necessary conditions for approximate solutions of vector as well as set optimization problems with variable domination structures.

On the one hand, we prove a new Ekeland's type variational principle suited for the problems under the study and we employed this variational principle to derive necessary conditions for approximate solutions of problems where the objective mapping is single-valued. Even if the statement (iii) in Theorem \ref{thm:bestK} is shown in a scalarized form, it is very useful to derive necessary conditions for approximate solutions. For further research, it could be interesting to improve the existence result concerning an exact solution of a perturbed scalarized problem in the third condition of the variational principle by a corresponding formulation in terms of nondominated solutions of a perturbed problem wrt variable domination structure. 

On the other hand, we prove a generalized openness result for a sum of set-valued mappings which allows us to give optimality conditions also for the case of set-valued objective mappings. 

It would be interesting to employ the necessary conditions for approximately nondominated solutions of vector-valued as well as set-valued optimization problems wrt variable domination structures derived in our paper for developing numerical procedures for solving these problems.}

\end{document}